\documentclass[12pt]{amsart}
\usepackage[T1]{fontenc}
\usepackage{lmodern}
\usepackage{microtype}
\usepackage[a4paper,margin=2.4cm]{geometry}
\usepackage{amsmath,amssymb,mathtools}

\usepackage{enumitem}
\usepackage[dvipsnames]{xcolor}
\usepackage[hidelinks]{hyperref}
\hypersetup{
    colorlinks=true,
    linkcolor=Blue,
    filecolor=magenta,      
    urlcolor=RoyalBlue,
    citecolor=Maroon
}

\theoremstyle{plain}
\newtheorem{theorem}{Theorem}[section]
\newtheorem{proposition}[theorem]{Proposition}
\newtheorem{lemma}[theorem]{Lemma}
\newtheorem{corollary}[theorem]{Corollary}
\newtheorem{conjecture}[theorem]{Conjecture}

\theoremstyle{definition}
\newtheorem{definition}[theorem]{Definition}
\newtheorem{notation}[theorem]{Notation}

\newtheorem{remark}[theorem]{Remark}

\numberwithin{equation}{theorem}

\usepackage{thm-restate}
\makeatletter
\newcommand{\rstlabel}[1]{%
  \ifthmt@thisistheone\label{#1}\else\tag{\ref*{#1}}\fi}
\makeatother

\newcommand{\kk}{\Bbbk}

\newcommand{\Vir}{\operatorname{Vir}}
\newcommand{\Sa}{\operatorname{S}}
\newcommand{\Sg}{\mathfrak S}
\newcommand{\UU}{\operatorname{U}}

\newcommand{\ad}{\operatorname{ad}}
\newcommand{\Ann}{\operatorname{Ann}}
\newcommand{\sgn}{\operatorname{sgn}}
\newcommand{\Alt}{\operatorname{Alt}}

\newcommand{\gr}{\operatorname{gr}}
\newcommand{\NN}{\mathbb{N}}

\newcommand{\ZZ}{\mathbb{Z}}

\newcommand{\Dx}{\operatorname{Dx}}
\newcommand{\Pprim}{\operatorname{P.Prim}}
\newcommand{\del}{\partial}

\newcommand{\mc}[1]{\mathcal{#1}}
\newcommand{\mf}[1]{\mathfrak{#1}}

\newcommand{\calA}{\mathcal A}
\newcommand{\calB}{\mathcal B}
\newcommand{\Der}{\operatorname{Der}}
\newcommand{\g}{\mathfrak{g}}
\newcommand{\geps}{\mathfrak{g}_\epsilon}
\newcommand{\St}{\mathrm{St}}
\newcommand{\End}{\mathrm{End}}
\newcommand{\length}[1]{\left\lvert #1\right\rvert}

\newcommand{\mult}{\mathrm{mult}}
\newcommand{\Spec}{\mathrm{Spec}}
\newcommand{\Prim}{\mathrm{Prim}}

\newcommand{\Span}{\mathrm{span}}
\newcommand{\wt}{\operatorname{wt}}
\newcommand{\ord}{\operatorname{ord}}

\DeclareMathOperator*{\iprod}{\overset{\rightarrow}{\prod}}

\title[Central quotients of Virasoro-type enveloping algebras]{Nonzero central quotients of the Virasoro enveloping algebra are simple}

\author{Tuan Anh Pham}
\address{T. A. Pham: School of Mathematics, The University of Edinburgh, James Clerk Maxwell Building, Peter Guthrie Tait Road, Edinburgh, 
 EH9 3FD, United Kingdom}
\email{tuan.pham@ed.ac.uk}

\author{Haijun Tan}
\address{H. Tan: School of Mathematics and Statistics, Northeast Normal University, Changchun, Jilin, 130024, P. R. China.}
\email{tanhj9999@163.com}


\begin{document}
\begin{abstract}
We prove a conjecture of Petukhov and Sierra that the central quotient $\mathrm{U}(\mathrm{Vir})/(z-c)$ of the universal enveloping algebra of the Virasoro algebra is simple for all $c \neq 0$. We extend the result to the $N=1$ super-Virasoro algebras.

Our proof reveals a new connection between enveloping algebras and polynomial identity algebras. Although the enveloping algebra itself does not satisfy a polynomial identity, the adjoint action on each filtered piece does satisfy a standard identity. In particular, the Witt algebra $W$, as an adjoint module over itself, satisfies the standard identity $\mathrm{St}_4$, whereas $\mathrm{Vir}$ fails by a nonzero multiple of the central charge. Propagating this gap through the PBW filtration of $\mathrm{U}(\mathrm{Vir})/(z-c)$, for every nonscalar element, we construct an alternating adjoint operator that lowers its order by exactly one.

As applications, we compute the annihilators of all representations with nonzero central character, prove a strengthened form of Etingof's conjecture on nonsplit extensions of $W$-modules, and reduce the classification of prime ideals of $\mathrm{U}(\mathrm{Vir})$ to those of $\mathrm{U}(W)$. 
\end{abstract}
\maketitle

\section{Introduction}

This paper answers a fundamental question concerning the universal enveloping algebra of the Virasoro algebra, one of the most important infinite-dimensional Lie algebras in mathematics and physics. Throughout, let $\kk$ be a field of characteristic $0$. 

The \emph{Witt algebra} $W = \Der \kk[t, t^{-1}]$ is the Lie algebra of algebraic vector fields on $\kk^\times$ and the \emph{Virasoro algebra} $\Vir = W \oplus \kk z$ is its unique (up to isomorphism) nontrivial one-dimensional central extension. Since their appearance in string theory and conformal field theory, these algebras have been studied extensively. Their representation theory is closely related to affine Kac--Moody algebras, $W$-algebras, and vertex operator algebras.

Despite being studied for more than 50 years, the ideal structure of the enveloping algebra $\UU(\Vir)$ remains mysterious. In particular, not much is known about the quotient $\calA_c := \UU(\Vir)/(z-c)$ where $c \neq 0$. Only recently, Petukhov and Sierra conjectured that $\calA_c$ is a simple ring \cite[Conjecture 7.2.1]{PetukhovSierra2022}, after their result that $\Sa(\Vir)/(z-c)$ is Poisson simple \cite[Corollary 3.3.5]{PetukhovSierra2022}. The conjecture thus tests whether the guiding principle that ideals of $\UU(\g)$ mirror Poisson ideals of $\Sa(\g)$ extends to the infinite-dimensional setting.
\begin{conjecture}[{\cite[Conjecture 7.2.1]{PetukhovSierra2022}}]\label{conjecture: Poisson simple}
    If $c \neq 0$, then $\UU(\Vir)/(z-c)$ is simple. 
\end{conjecture}
We prove Conjecture~\ref{conjecture: Poisson simple}, in the following stronger form. 
 \begin{restatable}{theorem}{simpleVir}\label{thm: simple Vir}
    For every field $\kk$ of characteristic $0$ and $c \in \kk^\times$, every nonzero subspace of $\calA_c = \UU(\Vir)/(z-c)$ that is stable under the adjoint action of $W$ contains $\kk$. 
    
    In particular, $\calA_c =  \UU(\Vir)/(z-c)$ is a simple $\kk$-algebra if and only if $c \in \kk^\times$.
\end{restatable}

More importantly, our proof establishes a surprising connection between enveloping algebras and polynomial identity algebras, which have historically been independent strands of noncommutative ring theory. Although in characteristic zero, $\UU(\g)$ satisfies no nontrivial polynomial identity except when $\g$ is abelian by \cite{Latyvsev63, Bahturin74}, the adjoint action of $\g$ on the PBW pieces $F_\ell \UU(\g)$ does satisfy a standard identity. We exploit the fact that $W$ satisfies a standard polynomial identity while $\Vir$ fails by a nonzero multiple of the central charge to construct operators that lower the order of nonzero elements by exactly one.

One difficulty in studying $\UU(\Vir)$ is that it is not Noetherian by Sierra and Walton \cite{SierraWalton14, SierraWalton16}, so many classical noncommutative ring theory tools do not apply. Nevertheless, substantial progress has been made. Iyudu and Sierra showed that ideals of $\UU(W)$ and $\UU(\Vir)/(z-c)$ are large in the sense that their quotients have finite Gelfand--Kirillov dimension in \cite{IyuduSierra20}. Then in \cite{SanchezSierra23}, León Sánchez and Sierra proved that $\Sa(W)$ satisfies the ascending chain condition on radical Poisson ideals using model theory. Petukhov and Sierra classified the prime and primitive Poisson ideals of $\Sa(W)$ and $\Sa(\Vir)$ in \cite{PetukhovSierra2022}. The first author also constructed a Dixmier map from Poisson primitive ideals of $\Sa(\mathrm{Vir})$ to primitive ideals of $\UU(\mathrm{Vir})$ in \cite{pham2025orbit}. However, these Poisson-theoretic results do not lift to $\UU(\Vir)$. The first author studied ideals of $\UU(\Vir)$ through homomorphic images of the so-called orbit homomorphisms \cite{PhamTimmins25, pham2026ideals}, but these homomorphisms all vanish on $z$. 

On the other hand, Biswal and Sierra proved simplicity at nonzero level for affine Kac--Moody algebras using large commuting families and an embedded affine $\widehat{\mf{sl}}_2$ \cite{BiswalSierra25}, neither tool is available for the Virasoro algebra. In the recent update of \cite{FeiginKanelKhoroshkin}, Feigin,  Kanel-Belov, and Khoroshkin reduced the simplicity of $\calA_c$ to a conjecture of Etingof on certain $W$-modules and their extensions. We prove the strongest form of Etingof's conjecture in Theorem~\ref{thm: super Etingof conj}.

We now explain the idea of our proof. Our approach is built on the important observation that although $\UU(\g)$ does not satisfy a polynomial identity, the adjoint action of $\g$ on the PBW piece $F_\ell \UU(\g)$ does satisfy a standard identity whose degree grows linearly in $\ell$. 
\begin{restatable}{theorem}{tensorpower}\label{thm: tensor power}
    Let $\rho: \g \to \End(V)$ be a representation of a Lie algebra $\g$ such that 
    \[\St_{2k}(\rho(x_1), \dots, \rho(x_{2k})) = 0 \qquad \text{for all}\qquad x_1, \dots, x_{2k} \in \g.\]
    Then for all $\ell \geq 1$, the $\g$-modules $V^{\otimes \ell}$ and $\Sa^\ell(V)$ satisfy $\St_{(2k-2)\ell+2}$. In particular, if the adjoint representation of $\g$ satisfies $\St_{2k}$, then $\Sa^\ell(\g)$ and $F_\ell \UU(\mathfrak g)$ satisfy $\St_{(2k-2)\ell+2}$.
\end{restatable}
Theorem~\ref{thm: tensor power} is in the spirit of Kostant's generalization~\cite{Kostant81} of the Amitsur--Levitzki theorem~\cite{AmitsurLevitzki50} to representations of Lie algebra, and Regev's result~\cite{Regev72} that the tensor product of two polynomial identity (PI) algebras is again PI. Our proof of Theorem~\ref{thm: tensor power} adapts Rosset's exterior algebra argument~\cite{Rosset76}. 

We then apply Razmyslov's result \cite{Razmyslov85} that $W$ as an adjoint $W$-module satisfies $\St_4$, while the Virasoro cocycle breaks this identity by a nonzero multiple of $c$. Let $\calA_c = \UU(\Vir)/(z-c)$ and we identify the associated graded of $\calA_c$ with respect to the PBW basis with $\Sa(W)$, i.e., $\gr \calA_c = \Sa(W)$. Writing $D_m= \ad(e_m)$, we define in Definition~\ref{def: alternating adjoint operator} the alternating adjoint operator $\St_n(D_q)  = \St_n(D_{q_1}, \dots, D_{q_n}) \in \End_\kk (\calA_c)$ for $q \in \ZZ^n$.

We then propagate these features to $F_{\ell} \calA_c$ to construct an order reduction algorithm. Here we say $w \in \calA_c$ has \emph{order $\ell$} if $w \in F_{\ell} \calA_c \setminus F_{\ell-1} \calA_c$. By applying $\St_{2\ell+2}(D_q)$ with carefully chosen ``base-$B$'' parameters $q \in \ZZ^{2\ell+2}$, we reduce the order of an element in $F_{\ell} \calA_c$ by exactly one. It is easy to see that applying Theorem~\ref{thm: degree reduction} repeatedly gives us Theorem~\ref{thm: simple Vir}.

\begin{restatable}{theorem}{degreduction}\label{thm: degree reduction}
    Let $c \neq 0$ and let $w \in \calA_c$ of order $d \geq 1$. Let $\mu \in \ZZ$ such that $x_\mu$ appears in $\gr(w)$, and choose $q, n, M$ as in Notation~\ref{not: degree reduction}. Then we have 
    \begin{equation}\rstlabel{eq: degree reduction}
        \gr \big( \St_{2d+2}(D_q)(w) \big) = C_d \, \St_{2d-2}(D_{q'}) \big( \partial_{x_\mu} \gr(w) \big), 
    \end{equation}
    which is nonzero, so $\St_{2d+2}(D_q)(w)$ has order exactly $d-1$.
\end{restatable}
We then extend Theorem~\ref{thm: simple Vir} to the $N=1$ super-Virasoro algebras, that is to the Neveu--Schwarz and Ramond algebras. The odd part introduces no additional difficulty, and the argument follows analogously.
\begin{theorem}[Theorems ~\ref{thm: NSR simple 1} and ~\ref{thm: NSR simple 2}]\label{thm: simple NSR}
    Let $\g$ be the Neveu--Schwarz or the Ramond algebra with central element $z$. If $c \in \kk^\times$, every nonzero subspace of $\UU(\g)/(z-c)$ that is stable under the (super) adjoint action of $\g$ contains $\kk$. 
    
    In particular, $\UU(\g)/(z-c)$ is a simple $\kk$-algebra if and only if $c \neq 0$.
\end{theorem}

We believe that our method generalizes to more Lie algebras, including $N=2$ super-Virasoro algebra, Virasoro-related Lie algebras, and central extensions of Lie algebras of vector fields on affine curves; see Section~\ref{subsec: Future work}.

Theorem~\ref{thm: simple Vir} has substantial consequences for the representation theory of $\Vir$.
\begin{restatable}{corollary}{anncor}\label{cor: ann}
    Let $N$ be a nonzero $\Vir$-module and $(z-c)N = 0$ for some $c \neq 0$. Then
    \begin{equation*}
        \Ann_{\UU(\Vir)} N = (z-c).
    \end{equation*}
\end{restatable}
\begin{proof}
    This follows from Theorem~\ref{thm: simple Vir} as $(z-c)$ is a maximal ideal of $\UU(\Vir)$.
\end{proof}
Corollary~\ref{cor: ann} is a stronger version of \cite[Theorem 1.4]{IyuduSierra20} and 
shows that every nonzero $\Vir$-module where $z$ acts by a nonzero scalar $c \neq 0$, including any yet to be constructed, is faithful over $\calA_c$. In particular, this applies to Verma modules and all their quotients, including $L(c,h)$ and the vacuum modules $V^c$; all restricted modules, i.e. smooth $V^c$-modules of \cite{MazorchukZhao14}, including the Whittaker modules; the weight modules $L(c,h) \otimes V_{\alpha,\beta}$ \cite{ChenGuoZhao13}; the non-weight, non-restricted modules $\Omega(\lambda,b)\otimes V$, with $\Omega(\lambda,b)$ free of rank one over $\UU(\kk e_0)$ \cite{TanZhao16}; the staggered modules \cite{KytolaRidout09}; and every object of the Virasoro tensor categories of \cite{CJORY21}, including logarithmic modules. 

Theorem~\ref{thm: simple Vir} also holds for $\UU(\Vir) \otimes_{\kk[z]} \kk(z)$. As a result, we reduce the classification problem of prime and primitive ideals of $\UU(\Vir)$ in \cite{pham2025orbit} to those of $\UU(W)$ in Theorem~\ref{thm: ideal UVir}.
\begin{corollary}\label{cor: simplicity Vir tensor}
    The $\kk(z)$-algebra $\UU(\Vir) \otimes_{\kk[z]} \kk(z)$ is simple, and every nonzero $\kk(z)$-subspace of it that is stable under the adjoint action of $W$ contains $\kk(z)$.
\end{corollary}
\begin{proof}
    Let $z \otimes 1$ be the central element of $\Vir_{\kk(z)} = \Vir \otimes_\kk \kk(z)$. Then $\UU(\Vir) \otimes_{\kk[z]} \kk(z) \cong \UU(\Vir_{\kk(z)})/(z \otimes 1 - 1\otimes z)$. The result follows from Theorem~\ref{thm: simple Vir} over $\kk(z)$.
\end{proof}

\textbf{Organization of the paper:} In Section~\ref{sec: standard}, we recall the PBW filtration on the enveloping algebra of the Virasoro algebra and $N=1$ super-Virasoro algebra. We then construct the alternating adjoint operator in Definition~\ref{def: alternating adjoint operator} and prove Theorem~\ref{thm: tensor power}. In Section~\ref{sec: vir}, we prove the results for the Virasoro algebra, including Theorem~\ref{thm: degree reduction} and Theorem~\ref{thm: simple Vir}. We extend the results to the Neveu--Schwarz and Ramond algebras in Section~\ref{sec: super}. In fact, we could have combined Sections~\ref{sec: vir} and \ref{sec: super}, but such presentation would require additional notation and care for the odd generators, which obscures the main ideas. Finally, we discuss the relation of our work to \cite{FeiginKanelKhoroshkin} and prove the strongest form of Etingof's conjecture in Section~\ref{subsec: FKBK work}. In Section~\ref{subsec: conseq}, we reduce the classification of prime and primitive ideals of $\UU(\Vir)$ to those of $\UU(W)$. We expect our method to generalize to more Lie (super)algebras, and formulate several conjectures in Section~\ref{subsec: Future work}.

\textbf{AI disclosure:} This paper was written by humans for humans. The initial project and the proofs of the main results, e.g. Theorem~\ref{thm: simple Vir} and Theorem~\ref{thm: simple NSR}, were developed completely by the authors. While we were revising the draft of the paper, interesting discussions with Plamen Koshlukov led us to reformulate our results in the language of PI algebras. As non-experts in PI theory, we used AI for initial feasibility assessment and literature suggestions, including Rosset's proof of the Amitsur--Levitzki theorem. We then formulated the question that later became Theorem~\ref{thm: tensor power}, inspired by further discussions with Plamen Koshlukov.

The paper was then rewritten significantly. AI tools were used limitedly in preparing the final manuscript: for language editing and for preliminary drafts of proofs of some routine lemmas, such as Lemma~\ref{lem: one block} and Lemma~\ref{lem:injective}, based on our approach. We checked these suggestions independently, completed the arguments, provided necessary details, and substantially restructured the proofs. All mathematical content was developed and verified by the authors and we take full responsibility for all claims and arguments in the paper.

\textbf{Acknowledgement: } Tuan Anh Pham thanks Plamen Koshlukov for various interesting discussions on PI theory and his particular interest in Theorem~\ref{thm: tensor power}. We heartfully thank Susan Sierra for bringing the recent update of \cite{FeiginKanelKhoroshkin} to our attention, valuable comments on this manuscript and stimulating conversations.

\section{Standard polynomial identities on the PBW  filtration}\label{sec: standard}

In this section, we recall the Virasoro, Neveu--Schwarz and Ramond algebras and the PBW filtration on their enveloping algebras. We then recall the definition of standard polynomials from PI theory, which is used to construct the alternating adjoint operator. The main result of this section is Theorem~\ref{thm: tensor power} that each PBW piece $F_\ell \UU(\g)$ does satisfy a standard identity whose degree grows linearly in $\ell$. 

All vector spaces and Lie (super)algebras are over a field $\kk$ of characteristic zero unless otherwise specified, and we abbreviate $\otimes_{\kk}$ to $\otimes$. If $\mc{L}$ is a Lie (super)algebra, we denote the universal enveloping algebra of $\mc{L}$ by $\UU(\mc{L})$ and the (super)symmetric algebra by $\Sa(\mc{L})$.

\subsection{PBW filtration of the enveloping algebra of \texorpdfstring{($N=1$ super-)}{super}Virasoro algebra}
\begin{definition}
Let $W = \Der \kk[t, t^{-1}] = \kk[t, t^{-1}]\partial$ be the \emph{Witt algebra}, where $\partial = \frac{d}{dt}$. It is also the Lie algebra of algebraic vector fields on $\kk^\times$. The \emph{Virasoro algebra} $\Vir = W \oplus \kk z$ is the unique (up to isomorphism) nontrivial one-dimensional central extension of $W$, where $z$ is central. Then $\{e_m := t^{m+1}\del \mid m \in \ZZ\} \cup \{z\}$ forms a basis for $\Vir$, with the Lie bracket 
\begin{equation}
    [e_m, e_n] = (n-m) e_{m+n} + \kappa(m) \, \delta_{m+n,0} \, z, \qquad \text{where} \qquad \kappa(m) = \frac{m^3-m}{12}.
\end{equation}
\end{definition}

The $N=1$ super-Virasoro algebras arise as the infinitesimal symmetry algebras of two-dimensional superconformal field theories and are obtained by adjoining a fermionic (odd) generator to the conformal generators. There are two natural choices of fermionic boundary condition: the Neveu–Schwarz algebra with half-integral modes and the Ramond algebra with integral modes. We introduce them uniformly as follows.
\begin{definition}
    Fix $\epsilon \in \{0, \frac{1}{2}\}$ and let $\ZZ_\epsilon = \ZZ + \epsilon$. The \emph{$N=1$ super-Virasoro algebra} $\geps$ is the Lie superalgebra with even basis $\{e_m\}_{m\in \ZZ} \cup\{z\}$ and odd basis $\{G_r\}_{r \in \ZZ_\epsilon}$. The Lie bracket of $\geps$ is given by $[z, \geps] = 0$, i.e., $z$ is central, and 
    \begin{equation}\label{svir-bracket}
    \begin{aligned}
        [e_m, e_n] &= (n-m)e_{m+n}+ \kappa(m) \, \delta_{m+n,0} \, z, \\
        [e_m, G_r] &= \left(r-\frac{1}{2}m\right)G_{m+r}, \\ 
        [G_r,G_s] &=-2e_{r+s}+ \frac{4r^2 - 1}{12} \delta_{r+s,0} \, z.
    \end{aligned}
    \end{equation}
    The Lie superalgebras $\mathfrak{g}_{1/2}$ and $\mathfrak{g}_0$ are also known as the \emph{Neveu--Schwarz} and \emph{Ramond algebras} respectively. 
\end{definition}
It is clear that the even part of $\geps$ is $\Vir$ and the odd part of $\geps$ is $\mathcal{G}:=\Span_{\kk}\{G_r\mid r\in \ZZ_\epsilon\}$. Thus $\geps = \Vir \oplus \, \mc{G}$ and the enveloping algebra $\UU(\geps)$ decomposes as a vector space
\[\UU(\geps)\cong\UU(\Vir)\otimes (\bigwedge \mathcal{G}), \qquad \text{where } \bigwedge \mathcal{G}:=\bigoplus_{t\in\NN}(\bigwedge ^t\mathcal{G}).\]
Unless otherwise specified, all results in the paper hold for both values of $\epsilon$. For $c \in \kk$, let $\calA_c : = \UU(\Vir)/(z-c)$ and $\calB_c:=\UU(\geps)/(z-c)$. Abusing notation, we denote the image in $\calA_c$ (resp. $\calB_c$) of an element of $\UU(\Vir)$ (resp. $\UU(\geps)$) by the same symbol. 
\begin{definition}\label{def:pbw-filtration}
For any Lie algebra $\mathfrak{g}$ and $d\in\NN$, let
\[
 F_d \UU(\g):=\Span \bigl\{v_1 v_2 \cdots v_k
 \mid 0 \le k \le d,\ v_j \in \g \bigr\}.
\]
It is easy to see that $\{F_d \UU(\g) \mid d \in \NN\}$ forms an exhaustive ascending filtration on $\UU(\g)$, which is also known as the \emph{PBW (standard) filtration} on $\UU(\g)$. For all $c \in \kk$, this induces the following filtrations on the quotients
\begin{align*}
    F_d \calA_c &:= \bigl (F_d \UU(\Vir) + (z-c) \bigr) \big/ (z-c), \qquad 
    F_d \calB_c &:= \bigl (F_d \UU(\geps) + (z-c) \bigr) \big/ (z-c),
\end{align*}
and note that the induced filtration on $\calA_0 = \UU(W)$ coincides with the PBW filtration on $\UU(W)$. We thus call $\{F_d \calA_c\}$ and $\{F_d \calB_c\}$ the \emph{PBW filtration} on $\calA_c$ and $\calB_c$. 
\end{definition}

Let $R$ be $\calA_c$ or $\calB_c$.
By the PBW theorem for Lie (super)algebras, the set
\begin{equation}\label{eq:basis-B}
\bigl\{ e_{m_1} \cdots e_{m_r} G_{n_1} \cdots G_{n_s} \mid  m_1 \le \cdots \le m_r \in \ZZ, \, n_1 < \cdots < n_s \in \ZZ_\epsilon \bigr\}
\end{equation}
forms a $\kk$-basis of $R$ and the associated graded ring of $R$ is isomorphic to
\[\gr R = \Sa(W \oplus \mathcal{G}) = \Sa(W) \otimes \bigwedge \mathcal{G}, \]
where $\mc{G} = 0$ if $R = \calA_c$.
\begin{definition}\label{def: order and weight}
    Let $R$ be $\calA_c$ or $\calB_c$. With respect to the standard filtration in Definition~\ref{def:pbw-filtration}, a nonzero $w \in R$ has \emph{order} $d$ if $w \in F_d R \setminus F_{d-1} R$. We denote the order of $w$ by $\ord(w)$. 

    For $d \in \NN$, let $\gr_d: F_d R \to F_d R / F_{d-1} R$ be the quotient map and let $\gr(w) := \gr_d(w)$ if $w$ has order $d$. 
    Denote $x_{m_1} \cdots x_{m_r} y_{n_1} \cdots y_{n_s}$ for the image of $ e_{m_1} \cdots e_{m_r} G_{n_1} \cdots G_{n_s}$ in $\gr R$. 

    Note that $R$ is also graded by \emph{weight}, the eigenvalue of $\ad(e_0)$, and we denote the weight of a homogeneous element $r$ by $\wt(r)$. In particular, $\calA_c$ is $\ZZ$-graded and $\calB_c$ is $\frac{1}{2}\ZZ$-graded respectively. Each monomial in the PBW basis \eqref{eq:basis-B} is homogeneous with weight 
    \[\wt(e_{m_1} \dots e_{m_r} G_{n_1} \cdots G_{n_s}) = m_1 + \dots + m_r + n_1 + \dots + n_s.\]
    This grading is sometimes referred to in the literature as degree, a term we reserve for the degree of polynomials in $ \kk[t_1, \dots, t_k]$.
\end{definition}

\subsection{Standard polynomials and Rosset's identity}\label{sec:alternants}
In this subsection, we recall the definition of standard polynomial from PI algebras and then use Rosset's exterior algebra argument from \cite{Rosset76} to prove Theorem~\ref{thm: tensor power}. We then define the alternating adjoint operator and look at its action on a monomial in $\calA_c$ or $\calB_c$.

\begin{definition}\label{def: standard polynomial}
    For every associative algebra $R$ and $x_1, \dots, x_n \in R$, let 
    \begin{equation}
        \St_n(x_1, \dots, x_n) = \sum_{\sigma \in \Sg_n} \sgn(\sigma) x_{\sigma(1)} \dots x_{\sigma(n)} \in R
    \end{equation}
    be the \emph{standard polynomial} as in \cite{AmitsurLevitzki50}. For every set $J = \{j_1 < \dots < j_k\} \subseteq \ZZ$, let 
    \[\St_J (x) = \St_k(x_{j_1}, \dots, x_{j_k})=\sum_{\sigma \in \Sg_k} \sgn(\sigma) x_{j_{\sigma(1)}} \dots x_{j_{\sigma(k)}}
    ,\]
    with the convention that $\St_\emptyset = 1$.

    We say that \emph{$R$ satisfies $\St_n$} if $\St_n(x_1, \dots, x_n) = 0$ for all $x_i \in R$.
\end{definition}
For a Lie algebra $\g$ and a representation $\rho: \g \to \End(V)$, asking the image of $\UU(\g)$ in $\End(V)$ satisfies $\St_n$ is too strong. Recall that a polynomial $f$ is a \emph{weak polynomial identity} of a pair $(R, \mc{L})$, where $R$ is an associative algebra and $\mc{L} \subseteq R$ is a Lie subalgebra of $R$, if $f$ vanishes under every substitution of its variables by elements of $\mc{L}$ \cite{Razmyslov94, Drensky21}.
\begin{definition}\label{def: std adjoint polynomial}
    We say that a representation $\rho \colon \g \to \End(V)$ \emph{satisfies $\St_n$ as a $\g$-module} if $\St_n$ is a weak polynomial identity of the pair $(\End(V), \rho(\g))$, that is,
    \[
        \St_n(\rho(x_1), \dots, \rho(x_n)) = 0 \qquad \text{for all} \qquad x_1, \dots, x_n \in \g.
    \] 
\end{definition}

The following lemma was the main ingredient in Rosset's proof \cite{Rosset76} of the Amitsur--Levitzki theorem \cite{AmitsurLevitzki50}.
\begin{lemma}\label{lem: Rosset}
   Let $R$ be an associative algebra and let $x_1, \dots, x_n \in R$. Let $\Lambda$ be the exterior algebra on $\xi_1, \dots, \xi_n$. Let $X = \sum_{j} \xi_j \otimes x_j \in \Lambda \otimes R$. Then, for all $0 \leq k \leq n$, we have 
    \begin{equation*}
        X^k = \sum_{J \subseteq [n], \length{J}=k} \xi_J \otimes \St_J (x)
    \end{equation*}
    where $[n]=\{1,\ldots, n\}$ and $\xi_J = \xi_{j_1}\dots \xi_{j_k}$ for $J = \{j_1 < \dots < j_k\}$.
\end{lemma}
\begin{proof}
    In $\Lambda \otimes R$, we have 
\[
  X^k = \sum_{i_1, \dots, i_k\in [n]} \xi_{i_1} \cdots \xi_{i_k} \otimes x_{i_1} \cdots x_{i_k}.
\]
As $\xi_i$ anticommute, every term with a repeated index vanishes. A tuple of distinct indices can be written uniquely as $(j_{\sigma(1)}, \dots, j_{\sigma(k)})$ with $J = \{j_1 < \dots < j_k\}$ and $\sigma \in \Sg_k$. Thus $\xi_{j_{\sigma(1)}} \cdots \xi_{j_{\sigma(k)}} = \sgn(\sigma)\, \xi_J$, then grouping the terms by $J$ gives the Lemma.
\end{proof}

\begin{definition}
    For $n \in \NN$, let $[n]$ denote the set $\{1, \dots, n\}$.
    An \emph{ordered set partition} of $[n]$ is a tuple $\mc{J} = (J_1, \dots, J_\ell)$ of disjoint, possibly empty, subsets of $[n]$ whose union is $[n]$. Let $\sgn(\mc{J})$ be the sign of the unique permutation that lists elements of $J_1$ in an increasing order then those of $J_2$ in an increasing order and so on. The same definition applies when $[n]$ is replaced by any finite set of integers.
\end{definition}
The following lemma is obvious and motivates a notation for ordered set partition $\mc{J}$.
\begin{lemma}\label{lem: easy}
    Let $n, k_1, \dots, k_\ell \in \NN$ and $n = k_1 + \dots + k_\ell$. Then 
    $\left\lfloor \frac{n}{2}\right\rfloor = \left\lfloor \frac{k_1}{2}\right\rfloor + \dots + \left\lfloor \frac{k_\ell}{2}\right\rfloor$
    if and only if at most one of the $k_i$ is odd. 
     \qed
\end{lemma}
\begin{notation}\label{not: h}
    For $\mc{J}$ an ordered set partition of $[n]$, denote $h = \left\lfloor \frac{n}{2} \right\rfloor$ and $h_a = \left\lfloor\frac{\length{J_a}}{2}  \right\rfloor$ for all $a = 1, \dots, \ell$.
\end{notation}

The following Proposition is a direct consequence of Rosset's argument and we include the proof for completeness.

\begin{proposition}\label{prop: Rosset split}
    Let $R$ be an associative algebra. For all $1 \leq i \leq n$ and $1 \leq a \leq \ell$, let $y_i^{(a)} \in R$ such that 
    \[ y_i^{(a)} y_j^{(b)} = y_j^{(b)} y_i^{(a)} \quad \text{for all }i, j \text{ and } a \neq b.\]
    Let $Y_i = \sum_{a} y_i^{(a)}$ and use Notation~\ref{not: h}. Then
    \begin{equation}\label{eq: Rosset split eq}
        \St_n(Y_1, \dots, Y_n) = \sum_{\mc{J}} \sgn(\mc{J}) \frac{h!}{h_1! \dots h_\ell!}  \St_{J_1} (y^{(1)}) \dots \St_{J_\ell}(y^{(\ell)}),
    \end{equation}
    where the sum runs over all ordered set partitions $\mc{J}=(J_1 \dots J_\ell)$ of $[n]$ with at most one block $J_a$ of odd size, and $\St_{J_a}(y^{(a)})$ is defined in Definition~\ref{def: standard polynomial}.
\end{proposition}
\begin{proof}
    Let $\Lambda$ be the exterior algebra generated by $\xi_1, \dots, \xi_n$. Define $X_a = \sum_{i} \xi_i \otimes y_i^{(a)} \in \Lambda \otimes R$ and let $X = \sum_a X_a$. For $a \neq b$, we have 
    \[X_a X_b = \sum_{i, j} \xi_i \xi_j y_i^{(a)} y_j^{(b)} = -  \sum_{i, j} \xi_j \xi_i y_j^{(b)} y_i^{(a)} = -X_b X_a.\]
    Hence $X^2 = \sum_a X_a^2$, and each $X_a^2$ commutes with every $X_b$. Since $n - 2h \in \{0, 1\}$, we have
    \begin{align*}
        X^n &= \Big( \sum_a X_a^2 \Big)^{h} \Big( \sum_b X_b \Big)^{n-2h} = \sum_{h_1 + \dots + h_\ell = h} \frac{h!}{h_1! \cdots h_\ell!} \, X_1^{2h_1} \cdots X_\ell^{2h_\ell} \Big( \sum_b X_b \Big)^{n-2h} \\
            &= \sum
               \frac{h!}{\lfloor k_1/2 \rfloor! \cdots \lfloor k_\ell/2 \rfloor!}\,
               X_1^{k_1} \cdots X_\ell^{k_\ell},
    \end{align*}
    where in the last step, the sum runs over all $(k_1, \dots, k_\ell) \in \NN^\ell$ with $k_1 + \dots + k_\ell = n$ and at most one $k_a$ odd (by Lemma~\ref{lem: easy}). Then expand each $X^{k_a}_a$ by Lemma~\ref{lem: Rosset} and note the product $\xi_{J_1} \dots \xi_{J_\ell}$ vanishes unless $\mc{J} = (J_1, \dots, J_\ell)$ is an ordered set partition of $[n]$. In that case, $\xi_{J_1} \dots \xi_{J_\ell} = \sgn(\mc{J}) \xi_1 \dots \xi_n$. The proposition thus follows from Lemma~\ref{lem: Rosset}.
\end{proof}

The following theorem shows that standard identities propagate \emph{linearly} from a representation to its tensor powers. The result is in the spirit of Kostant's generalization ~\cite{Kostant81} of the Amitsur--Levitzki theorem for Lie algebras. For any vector space $V$, let $\Sa^\ell(V)$ be the $\ell$th symmetric power of $V$. 

\tensorpower*
\begin{proof}
    It is easy to see that $V$ also satisfies $\St_{k'}$ for every $k' \geq 2k$. Let $x_i^{(a)}$ denote the action of $x_i$ on the $a$-th tensor factor. 
    By Proposition~\ref{prop: Rosset split} and Lemma~\ref{lem: easy}, all $J_a$ appearing in~\eqref{eq: Rosset split eq} have even length. Thus, there is a block $J_{b}$ with $\length{J_{b}} \geq 2k$, so $\St_{J_{b}}(x^{(b)}) = 0$. Then by Proposition~\ref{prop: Rosset split}, $V^{\otimes \ell}$ satisfies $\St_{(2k-2)\ell+2}$ and so does $\Sa^\ell(V)$.

    By assumption, $F_1\UU(\g) = \kk \oplus \g$ satisfies $\St_{2k}$. As the $\g$-module homomorphism $(F_1\UU(\g))^{\otimes \ell} \rightarrow F_\ell \UU(\g)$ is surjective, $F_\ell \UU(\g)$ satisfies $\St_{(2k-2)\ell+2}$, and so does $\Sa^\ell(\g)$.
\end{proof}
Theorem~\ref{thm: tensor power} is an analog, in the weaker setting of Definition~\ref{def: std adjoint polynomial}, of Regev's theorem~\cite{Regev72} that tensor products of PI algebras are PI. Whereas \cite{Regev72} only shows the existence of some polynomial identity, we show that the diagonal action of $\g$ on $V^{\otimes \ell}$ satisfies an explicit standard identity whose degree is \emph{linear} in~$\ell$, and this bound is sharp (for $k = 2$; see Remark~\ref{rem: bound is tight}).

The following lemma generalizes an elementary property of primitive elements
$x$ in a Hopf algebra, namely $\Delta^{(\ell)}(x) = \sum_{a} x^{(a)}$.

\begin{lemma}\label{lem: mult Leibniz}
    Let $V$ be a subspace of an associative algebra $R$ and $\delta_1, \dots, \delta_n$ be derivations of $R$ such that $\delta_i(V) \subseteq V$. Let $\mult^\ell: V^{\otimes \ell} \to R$ be the multiplication map. Let $\delta_i^{(a)} \in \End(V^{\otimes \ell})$ denote $\delta_i|_V$ acting on the $a$-th tensor factor, and let $Y_i = \sum_{a=1}^\ell \delta_i^{(a)}$. Then $\delta_i^{(a)}$ commutes with $\delta_{i'}^{(b)}$ for all $a \neq b$, and  
    \begin{equation}\label{eq: mult Leibinz}
        \St_n(\delta_1, \dots, \delta_n) \circ \mult_\ell = \mult_\ell \circ \St_n(Y_1, \dots, Y_n).
    \end{equation}
\end{lemma}
\begin{proof}
    It is clear that operators acting on different tensor factors commute with each other. By Leibniz rule, $\delta_i \circ \mult_\ell = \mult_\ell \circ Y_i$ so $\delta_{i_1} \dots \delta_{i_k} \circ \mult_\ell = \mult_\ell \circ Y_{i_1} \dots Y_{i_k}$ by induction. Thus \eqref{eq: mult Leibinz} follows from taking alternating sum.
\end{proof}

We now apply these results to the adjoint action of $\Vir$ and $\geps$ on $\calA_c$ and $\calB_c$ respectively.

\begin{definition}\label{def: alternating adjoint operator}
Let $c \in \kk$ and let $R$ be $\calA_c$ or $\calB_c$. For $m \in \ZZ$, let $D_{m} = \ad(e_{m}) \in \End_\kk(R)$, which is an even derivation of $R$. For $q = (q_1, \dots, q_n) \in \ZZ^n$ and $J = \{j_1 < \dots < j_k\} \subseteq [n]$, let $q_J = (q_{j_1}, \dots, q_{j_k})$ and 
\[
\St_J(D_q) = \St_k(D_{q_{j_1}}, \dots, D_{q_{j_k}}).
\]
We call the linear map $\St_J(D_q): R \to R$ the \emph{alternating adjoint operator} corresponding to $q_J$, and write $\St_n(D_q) = \St_{[n]}(D_q)$. 
\end{definition}

Since $\ad$ extends to an algebra map $\UU(\g) \to \End_\kk(R)$ and $D_{q_1} D_{q_2} = \ad(e_{q_1} e_{q_2})$, we have
\begin{equation}\label{eq: St ad}
  \St_J(D_q) = \ad\big(\St_J(e_q)\big).
\end{equation}
As $\St_J(D_q)$ maps homogeneous elements of weight $m$ to homogeneous elements of weight $m+ \sum_{j \in J} q_j$, we also denote $\wt(J) = \sum_{j \in J} q_j$.

\begin{corollary}\label{cor: apply adjoint}
    Let $c \in \kk$ and $R$ be $\calA_c$ or $\calB_c$. Let $q = (q_1, \dots, q_n) \in \ZZ^n$ where $n=2h$ is even. For all $v_1, \dots, v_\ell \in F_1 R$, we have 
    \[\St_n(D_q) (v_1 \dots v_\ell) = \sum_{\mc{J}} \sgn(\mc{J}) \frac{h!}{h_1! \dots h_\ell!}  \mult _\ell 
    \Bigl(\bigotimes_{a=1}^{\ell}\St_{J_a}(D_q)(v_a)\Bigr),\]
    where the sum runs over all ordered set partitions $\mc{J}=(J_1 \dots J_\ell)$ of $[n]$ with $\length{J_a} = 2h_a$ even for all $a$. 
\end{corollary}
\begin{proof}
    Since $D_m(\kk) = 0$ and 
    \begin{align}
        D_{m}(e_p) = (p-m)e_{p+m} +  \kappa(m) \, \delta_{p+m, 0} \, c, \qquad D_m(G_r) = \Bigl(r - \frac{m}{2}\Bigr) G_{r+m},\label{eq: ad 1}
    \end{align}
    it is easy to see that $V = F_1 R$ is stable under every $\delta_i = D_{q_i}$, which allows us to apply Lemma~\ref{lem: mult Leibniz}. By Proposition~\ref{prop: Rosset split} and Lemma~\ref{lem: easy}, all $J_a$ appearing have even length.
\end{proof}

\section{Central quotients of the Virasoro enveloping algebra}\label{sec: vir}
In this section, we will prove Theorem~\ref{thm: simple Vir} that $\calA_c$ is a simple ring for all $c \neq 0$. It suffices to show that every element $w \in \calA_c$ of order $d$ can be reduced to a nonzero scalar. The main idea is to leverage the gap between the Witt algebra $W$, which satisfies the standard Lie identity of degree $5$, and the Virasoro algebra $\Vir$, which fails by a multiple of the central term. For $c \neq 0$, this central term lets us construct an alternating adjoint operator $\St_n(D_q)$ that lowers the order of $w$ by exactly one and in a controlled way.

We first apply the alternating adjoint operator $\St_n(D_q)$ in Definition~\ref{def: alternating adjoint operator}, with $n$ even, to a monomial $e_{m_1} \cdots e_{m_\ell}$. Corollary~\ref{cor: apply adjoint} implies that these actions are distributed among the letters and each letter receives an even number of actions. Thus, each letter $e_p$ can absorb $0$ actions (remains unchanged), $2$ actions (shifted to another basis element $e_{p+ \wt(J)}$, possibly plus a scalar) or $4$ actions (deleted and replaced by a scalar). This is proven in Proposition~\ref{prop: monomial applied}.

Then we choose $n = 2\ell+2$. In this case, at least one letter must absorb $4$ actions, so it is replaced by a scalar, which lowers the order by $1$. Finally, we choose $q \in \ZZ^n$ so that exactly one letter, a fixed $e_\mu$, is replaced by a nonzero scalar and $\St_n(D_q)$ kills every monomial of order $<\ell$. The $q_i$ are also chosen as different powers of a large base $B$ and they determine the new index of every remaining letter so distinct terms cannot cancel and surviving coefficients are nonzero.
\subsection{The action of the alternating adjoint operator on a letter}
In this subsection, we show the dichotomy of the action of alternating adjoint on an element $e_p \in F_1 \calA_c$. We recover the fact that $W$ as an adjoint module over itself satisfies $\St_4$ and prove that $\Vir$ fails $\St_4$ as a $W$-module by a nonzero scalar. 

\begin{definition}  
    Let $T = (t_1, \dots, t_k)$ be commuting indeterminates and we define the \emph{Vandermonde determinant} to be 
    \begin{equation}\label{eq: vandermonde}
    \Delta (t_1, \dots, t_k) := \prod_{1 \leq a < b \leq k} (t_b-t_a).
    \end{equation}
    For all $f \in \kk[t_1, \dots, t_k]$, we denote 
    \[\Alt_k f = \sum_{\sigma \in \Sg_k} \sgn(\sigma) f(t_{\sigma(1)}, \dots, t_{\sigma(k)})\]
    A polynomial  $f \in \kk[t_1, \dots, t_k]$ is called \emph{alternating} if for all $\sigma \in \Sg_k$, we have 
    \[f(t_{\sigma(1)}, \dots, t_{\sigma(k)}) = \sgn(\sigma) f(t_1, \dots, t_k).\]
\end{definition}

It is well-known that the Vandermonde determinant divides every alternating polynomial in $t_1, \dots, t_k$.
\begin{lemma}\label{lem: div vandermonde}
    If $f \in \kk[t_1, \dots, t_k]$ is an alternating polynomial, then $\Delta (t_1, \dots, t_k) | f$. In particular, if $\deg f < \binom{k}{2}$, then $f = 0$.
    \qed
\end{lemma}

We now define two polynomials that will be used repeatedly. The lemma below says that they appear naturally as coefficients of successive application of the adjoint action on the generators $e_p$.
\begin{definition}
    For $J =\{j_1 < \dots < j_k\}\subseteq [n]$ and $q \in \ZZ^n$, we denote $\Delta(q_J) = \Delta(q_{j_1}, \dots,q_{j_k})$.  For $p \in \ZZ$, let 
    \begin{equation}
        A_k(p; T) = \prod_{\alpha=1}^k (p - t_\alpha + \sum_{\beta > \alpha} t_\beta), \quad B_k(p;T) = \kappa(t_1)A_{k-1}(p; t_2, \dots, t_k).
    \end{equation}
\end{definition}

\begin{lemma}\label{lem: term 1 ele}
    For all $p \in \ZZ$ and $\nu=(\nu_1, \dots, \nu_k) \in \ZZ^k$, we have 
    \begin{equation}
        D_{\nu_1} \dots D_{\nu_k} (e_p) = A_k(p; \nu) \, e_{p + \wt(\nu)} +  B_{k}(p; \nu) \, \delta_{p+ \wt(\nu), 0} \, c.  
    \end{equation}
\end{lemma}
\begin{proof}
    We prove the lemma by induction. The case $k=1$ follows from \eqref{eq: ad 1}. Let $p' = p+ \nu_2 + \dots + \nu_k$. Assuming the induction hypothesis, we have that $D_{\nu_2} \dots D_{\nu_k} (e_p) = A_{k-1} (p; \nu_2, \dots, \nu_k) e_{p'} + \gamma$ for some $\gamma \in \kk$. Since $D_{\nu_1}(\kk) = 0$, from \eqref{eq: ad 1}, we have 
    \[D_{\nu_1} \dots D_{\nu_k} (e_p) = A_{k-1} (p; \nu_2, \dots, \nu_k) \Big( (p' - \nu_1) e_{p' + \nu_1} + \kappa(\nu_1) \delta_{p' + \nu_1, 0} \, c \Big). \]
    It is easy to check that $A_{k-1} (p; \nu_2, \dots, \nu_k)(p' - \nu_1) = A_k (p; \nu)$ and $\kappa(\nu_1) A_{k-1} (p; \nu_2, \dots, \nu_k) = B_k (p; \nu)$, the lemma thus follows. 
\end{proof}

We now show the main result of this subsection: the action of the standard polynomial on basis elements of $\Vir$. Most of the statements are well known; the main purpose of the lemma is to establish precise control of the scalar $\kappa_J (p)$ in Lemma~\ref{lem: one block}(5).
\begin{lemma}\label{lem: one block}
    Let $J = \{j_1 < \dots < j_k\} \subseteq [n]$ and $p \in \ZZ$. Then 
    \begin{equation}\label{eq: std ad 1}
        \St_J(D_q)(e_p) = \omega_J (p) e_{p+\wt(J)} +  \kappa_J (p) \, c,
    \end{equation}
    where
    \[\omega_J(p) = \Alt_k A_k(p; q_J), \qquad \kappa_J(p) = \delta_{p+wt(J), 0} \, \Alt_k B_k(p; q_J).\]    
    Moreover, we have 
    \begin{enumerate}[label=\textup{(\arabic*)}]
        \item if $k \geq 4$, then $\omega_J (p) = 0$;
        \item if $J = \{i < i'\}$, then $\omega_J(p) =(q_{i'} - q_i) (p - q_i - q_{i'})$;
        \item if $p+\wt(J) \neq 0$, then $\kappa_J(p) = 0$;
        \item if $k \ge 5$, then $\kappa_J(p) = 0$; 
        \item if $k = 4$ and $p+\wt(J) = 0$, then $\kappa_J(p) = \frac{1}{6}\Delta(q_J)$.
    \end{enumerate}
\end{lemma}
\begin{proof}
    Equation \eqref{eq: std ad 1} follows easily from Lemma~\ref{lem: term 1 ele} and (3) is obvious. As $\Alt_k A_k(p; T)$ is a polynomial of total degree at most $k$, $\omega_J (p) = 0$ for all $k \geq 4$  by Lemma~\ref{lem: div vandermonde} and we have (1). Similarly, $\Alt_k B_k(p; T)$ is a polynomial of total degree at most $k+2$, so $\kappa_J (p) = 0$ for all $k \geq 5$ by Lemma~\ref{lem: div vandermonde}, giving (4). For (2), a direct calculation shows that
    \[A_2(p; q_i, q_{i'}) - A_2(p; q_{i'}, q_{i}) = (q_{i'} - q_i) (p - q_i - q_{i'}).\]
    We only need to show (5). Substituting $p = -(t_1 + \dots + t_4)$, gives $\kappa_J(p) = (\Alt_4 g)(q_J)$ where 
    \[
        g(t) = -\kappa(t_1)\,(t_1 + 2t_2)(t_1 + t_2 + 2t_3)(t_1 + t_2 + t_3 + 2t_4).
    \]   
    Recall that $\kappa(t_1) = \frac{1}{12} t_1^3 - \frac{1}{12} t_1$. The term $\frac{1}{12} t_1$ contributes an alternating polynomial of degree $4$, which is zero by Lemma~\ref{lem: div vandermonde}. On the other hand, the term $\frac{1}{12}t_1^3$ gives a homogeneous polynomial $g_6$ of degree $6$, so together, we have
    \[\Alt_4 g_6 = \gamma\, \Delta(t_1, t_2, t_3, t_4) \qquad \text{for some } \gamma \in \kk.\]
    It suffices to compare the coefficient of $t_1^3 t_2^2 t_3$ of both sides. The coefficient of $t_1^3 t_2^2 t_3$ in $\Delta(t)$ is $1$; while in $\Alt_4 g_6$, only three permutations contribute: the identity, the transposition swapping $t_2$ and $t_3$, and the transposition swapping $t_3$ and $t_4$. So $\gamma = -\frac{1}{12}(6 - 4 - 4) = \frac{1}{6}$.
\end{proof}

\begin{remark}\label{rem: std W Vir}
    Lemma~\ref{lem: one block} implies that $\St_4(\ad(e_{q_1}), \dots, \ad(e_{q_4})) = 0$ in $W$, which recovers the well-known result that $W$ satisfies the standard Lie identity of degree 5 \cite{Razmyslov85}. The lemma also shows that in $\Vir$ the identity fails with only the central scalar $c \kappa_J$.
\end{remark}

\subsection{Order reduction algorithm}
In this subsection, we propagate the features in Remark~\ref{rem: std W Vir} to a monomial in Proposition~\ref{prop: leading term}. The main ingredient is the choice of parameters of $\St_{2\ell+2}(D_q)$ in Notation~\ref{not: degree reduction} that simplifies bookkeeping and ensures the result is nonzero (Proposition~\ref{prop: leading term}). As a consequence, we obtain the order reduction algorithm in Theorem~\ref{thm: degree reduction} and simplicity of $\calA_c$ in Theorem~\ref{thm: simple Vir}. 

Since we are working over the noncommutative ring $\calA_c$, the order of multiplication matters, so we introduce the following notation.
\begin{notation}
    For all $A = \{a_1 < \dots < a_k\}$ and elements $v_a$ in an associative algebra $R$, let 
    \[\iprod_{a \in A} v_a := v_{a_1} v_{a_2} \dots v_{a_k}\] denote the product taken in increasing order of the index.
\end{notation}

We first show a general vanishing statement. This follows from Corollary~\ref{cor: apply adjoint} and easy book-keeping.

\begin{proposition}\label{prop: monomial applied}
    Let $c \in \kk$ and let $q = (q_1, \dots, q_n) \in \ZZ^n$ for $n = 2h$ even. Then for all $u = e_{m_1} \dots e_{m_\ell} \in \calA_c$, we have 
    \begin{equation}\label{eq: std product}
        \St_n(D_q)(u) = \sum_{\mc{J} \in \mc{P}(u)} \sgn(\mc{J})\,
      \frac{h!}{h_1! \dots h_\ell!} \, \iprod_{a \in [\ell]} \St_{J_a}(D_q)(e_{m_a}),
    \end{equation}
    where $h_a = \length{J_a}/2$ and $\mc{P}(u)$ is the set of ordered set partitions $\mc{J} = (J_1, \dots, J_\ell)$ of $[n]$ such that $\length{J_a} \in \{0, 2, 4\}$ and $m_a + \wt(J_a) = 0$ if $\length{J_a} = 4$ for all $a$. Moreover, 
  \[
    \St_{J_a}(D_q)(e_{m_a}) =
    \begin{cases}
      e_{m_a}, & \text{if } J_a = \emptyset; \\
      \omega_{J_a}(m_a)\, e_{m_a + \wt(J_a)} +  \kappa_{J_a}(m_a)\, c,
        & \text{if } \length{J_a} = 2; \\
      \kappa_{J_a}(m_a)\, c, & \text{if } \length{J_a} = 4.
    \end{cases}
  \]
\end{proposition}
\begin{proof}
    Corollary~\ref{cor: apply adjoint} gives the initial formula for Equation \eqref{eq: std product}. Then Lemma~\ref{lem: one block} (1) and (4) show that for blocks $J_a$ with $\length{J_a} \ge 5$, we have $\St_{J_a}(D_q)(e_{m_a}) = 0$. Thus, the sum in \eqref{eq: std product} runs over all ordered set partitions $\mathcal J = (J_1, \dots, J_\ell)$ of $[n]$ such that $\length{J_a} \in \{0, 2, 4\}$. Furthermore, if $\length{J_a} =4$ and $m_a + \wt(J_a) \neq 0$, by Lemma~\ref{lem: one block} (1) and (3), we have $\St_{J_a}(D_q)(e_{m_a}) = 0$.  
\end{proof}

\begin{remark}\label{remark: bound 1}
    If $c=0$, blocks of size $4$ also vanish, so $\St_n(D_q)(u) = 0$ if $n > 2\ell$. This is Theorem~\ref{thm: tensor power} for $\g = W$. 
\end{remark}

Thus, we need to control the blocks $J_a$ of size $4$ and the choice of $q$ below ensures that at most one block of size $4$ remains at a particular position.

\begin{notation}\label{not: degree reduction}
    Let $c \neq 0$ and let $w \in \calA_c$ or $\calB_c$ of order $d \geq 1$. Let $M$ be an upper bound for the absolute values of all indices appearing in $w$ written in the PBW basis. Let $n = 2d+2$ and $B \in \NN$ such that $B > 2M+1$. Let $\mu \in \ZZ$ such that $\length{\mu} \leq M$ and define 
    \begin{equation}
        q_i =  \begin{cases}
            B^i & \text{ if } 1 \leq i \leq n, i \neq 4; \\ 
            -(q_1+q_2+q_3)-\mu & \text{ if } i = 4.
        \end{cases}
    \end{equation}
    Thus $q_1+\dots + q_4 = -\mu$ so the block $\{1, 2, 3, 4\}$ only acts on $e_\mu$. Let 
    \begin{equation}\label{eqdef: C}
        q' = (q_5, \dots, q_n) \qquad \text{and} \qquad C_d = \frac{c}{12} d(d+1) \Delta (q_1, \dots, q_4).
    \end{equation}
    For a monomial $u = e_{m_1} \dots e_{m_\ell} \in R$ and $1 \leq a \leq \ell$, let 
    \[u_{\hat{a}} = e_{m_1} \dots \widehat{e_{m_a}} \dots e_{m_\ell}.\]
\end{notation}

The following lemma is stated for indices in $\frac{1}{2}\ZZ$ in order to be reused for the $N=1$ super-Virasoro algebra in Section~\ref{sec: super}. It also explains the choices of $n$ and $q_i$ in Notation~\ref{not: degree reduction}.

\begin{lemma}\label{lem: conseq of choice}
    In Notation~\ref{not: degree reduction}: 
    \begin{enumerate}
        \item Let $m \in \ZZ$ with $\length{m} \leq M$ and let $\emptyset \neq J \subseteq [n]$. Then 
        \[m + \wt(J) =0 \Leftrightarrow m = \mu, J = \{1, 2, 3, 4\}.\]
        \item Let $m, m' \in \frac{1}{2}\ZZ$ such that  $\length{m}, \length{m'} \leq M$ and $m-m' \in \ZZ$. Let $J, J' \subseteq \{5, \dots, n\}$. 
        If $m+\wt(J) = m'+\wt(J')$, then $J = J'$ and $m = m'$.
        \item For every 2-element set $J \subseteq \{5, \dots, n\}$, we have $\wt(J) > 2M$ and $\omega_J(m) \neq 0$ for every $m \in \ZZ$ with $\length{m} \leq M$. Moreover $\Delta(q_1, \dots, q_4) \neq 0$. 
    \end{enumerate}
\end{lemma}
\begin{proof}
    Note that $\wt(\{1, 2, 3, 4\}) = q_1 + q_2 + q_3 + q_4 = -\mu$, so $m + \wt(\{1, 2, 3, 4\}) = m - \mu$ vanishes if and only if $m = \mu$. Let $J \neq \{1, 2, 3, 4\}$. If $4 \notin J$, then $\wt(J) \geq B > \length{m}$ and $m + \wt(J) > 0$. Otherwise if $4 \in J$, let $J' = J \setminus \{4\}$. Then $m + \wt(J) = (m - \mu) + S$ where 
    \begin{align*}
        S = \sum_{j \in J'} B^j - (B + B^2 + B^3).
    \end{align*}
    It is easy to see that $S$ is a nonzero multiple of $B$ and $|S| \ge B > \length{m-\mu}$. Hence $m + \wt(J) \neq 0$ and we have (1).
    
    For (2), the assumption implies $\wt(J) - \wt(J') = m' - m$. Since $\wt(J) - \wt(J')$ is a multiple of $B$ and $\length{m'-m} \leq 2M < B$, both sides are zero. Thus $m = m'$, $\wt(J) = \wt(J')$. It easily follows that $J = J'$.
    
    Finally, for (3), Lemma~\ref{lem: one block}(2) implies that 
    \[
        \omega_J(m) = (q_{i'} - q_i)\bigl(m - \wt(J) \bigr).
    \]
    It is easy to see that $\omega_J(m) \neq 0$ since $m \leq M \leq 2M < \wt(J)$. As $0< q_1 < q_2 < q_3$ and $q_4 = -(B + B^2 + B^3) - \mu < 0$, the Vandermonde determinant $\Delta(q_1, \dots, q_4)$ is nonzero.
\end{proof}

The following Proposition deletes the letter $e_\mu$ appearing in $u$, where $\length{\mu} \leq M$. 

\begin{proposition}\label{prop: leading term}
    Let $u = e_{m_1} \dots e_{m_\ell}$ be an ordered monomial with $\length{m_a} \leq M$ for all $a$, and use Notation~\ref{not: degree reduction}. Then $\St_{2d+2}(D_q)(u) = 0$ if $\ell <d$ and 
    \begin{equation}
            \St_{2d+2}(D_q)(u) = C_d \sum_{a: m_a = \mu} \St_{2d-2}(D_{q'})(u_{\hat{a}}) \qquad \text{if}\qquad \ell = d.
        \end{equation}    
    In particular, $\St_{2d+2}(D_q)(u) \in F_{d-1} \calA_c$ and $\St_{2d+2}(D_q)(u) = 0$ if there exists no $a$ such that $m_a = \mu$.  
\end{proposition}
\begin{proof}
    Let $\mc{J} \in \mc{P}(u)$ as in Proposition~\ref{prop: monomial applied}. If $\length{J_a} = 4$ then $m_a + \wt(J_a) =0$ so $J_a = \{1, 2, 3, 4\}$ and $m_a = \mu$ by Lemma \ref{lem: conseq of choice}(1). Hence at most one block has size $4$ and $2d+2 = \sum \length{J_a} \leq 2\ell+2$. So $\St_{2d+2}(D_q)(u) = 0$ if $\ell < d$.

    For $\ell = d$, equality forces $J_{a_0} = \{1, 2, 3, 4\}$ for some $a_0$ with $m_{a_0} = \mu$ and $\widehat{\mc{J}} = (J_a)_{a \neq a_0}$ is an ordered partition of $\{5, \dots, n\}$ into pairs. Note that $\sgn(\widehat{\mc{J}}) = \sgn(\mc{J})$, so by Proposition~\ref{prop: monomial applied}, $\St_{2d+2}(D_q)(u)$ equals
    \begin{equation*}
         \frac{(d+1)!}{2}\frac{c}{6} \Delta(q_1, \dots, q_4) \sum_{a_0: m_{a_0} = \mu} \Big(\sum_{\widehat{\mc{J}}} \sgn(\widehat{\mc{J}}) \prod_{a \neq a_0} \omega_{J_a}(m_a) \iprod_{a \neq a_0} e_{m_a + \wt(J_a)}  \Big).
    \end{equation*}
    
    Apply Proposition~\ref{prop: monomial applied} to $\St_{2d-2}(D_{q'})$ with index set $\{5, \dots, n\}$ and $u_{\hat{a}_0}$. Note that $u_{\hat{a}_0}$ has order $d-1$, so $\{5, \dots, n\}$ has to split into $d-1$ blocks of size at most $2$. So similarly, 
    \begin{equation*}
        \St_{2d-2}(D_{q'}) \cdot 
        u_{\hat{a}_0} = (d-1)!\sum_{\widehat{\mc{J}}}\sgn(\widehat{\mc{J}}) \prod_{a \neq a_0} \omega_{J_a}(m_a) \iprod_{a \neq a_0} e_{m_a + \wt(J_a)},
    \end{equation*}
    and we obtain the Proposition.
\end{proof}
It remains to show that $\St_{2d+2}(D_q)(w) \neq 0$ when $w$ has order $d$, which will follow from an injectivity property of $\St_{2d-2}(D_{q'})$. Since $D_m$ preserves the PBW filtration and the central term of $D_m(e_p)$ lies in $F_0 \calA_c$, it induces the Poisson derivation $\{x_m, \_\}$ on $\gr \calA_c \cong \Sa(W)$ defined by $x_p \mapsto (p-m)x_{p+m}$. Thus $\St_{2d-2}(D_{q'})$ induces a linear map 
\[\St_{2d-2}(D_{q'}): \Sa^{d-1}(W) \to \Sa^{d-1}(W).\]

\begin{lemma}\label{lem:injective}
    Use Notation ~\ref{not: degree reduction} and let 
    \[\mc{X}_{M,d} = \Span\{x_{k_1}\dots x_{k_{d}} \in \Sa^{d}(W) \text{ such that } \length{k_a} \leq M \}.\]
    Then the map $\St_{2d-2}(D_{q'})$ is injective on $\mc{X}_{M, d-1}$.
\end{lemma}
\begin{proof}
    Let $X  = x_{k_1}\dots x_{k_{d-1}} \in \mc{X}_{M, d-1}$. The proof of Proposition~\ref{prop: leading term} with $u_{\hat{a}_0}$ replaced by $x_{k_1}\dots x_{k_{d-1}}$ gives
    \begin{equation}\label{eq: injective help 1}
        \St_{2d-2}(D_{q'})(X) = (d-1)! \sum_{\widehat{\mc{J}}} \sgn(\widehat{\mc{J}}) \prod_{a=1}^{d-1}\omega_{J_a}(k_a) \prod_{a=1}^{d-1} x_{k_a + \wt(J_a)},
    \end{equation}
    where $\widehat{\mc{J}}$ runs over the ordered partitions of $\{5, \dots, n\}$ into pairs. Lemma~\ref{lem: conseq of choice}(2) implies that $k_a + \wt(J_a)$ determines both $k_a$ and $J_a$. So each monomial $Y$ appearing in $\St_{2d-2}(D_{q'})(X)$ determines $X$. Therefore the images of distinct monomials in $\mc X_{M,d-1}$ have disjoint supports.
    
    Moreover, suppose $\widehat{\mc{I}}$ and $\widehat{\mc{J}}$ give the same monomial $Y$ in  $\St_{2d-2}(D_{q'})(X)$. Then by the previous paragraph, $\widehat{\mc{I}}$ is obtained from $\widehat{\mc{J}}$ by permuting blocks $J_a$ with the same $k_a$. Such a permutation does not change $\prod_{a} \omega_{J_a}(k_a)$ or $\sgn(\widehat{\mc{J}})$. So $\St_{2d-2}(D_{q'})(X) \neq 0$, by Lemma~\ref{lem: conseq of choice}(3). Thus $\St_{2d-2}(D_{q'})$ is injective on $\mc{X}_{M, d-1}$.
\end{proof}

We now prove the order reduction algorithm. It follows mainly from Proposition~\ref{prop: leading term} and Lemma~\ref{lem:injective}. We denote the operator $\frac{\partial}{\partial {x_\mu}}$ on $\Sa(W) = \kk[x_m \mid m \in \ZZ]$ by $\partial_{x_\mu}$.

\degreduction*
\begin{proof}
    Write $w = \sum_{\lambda} \gamma_\lambda u_\lambda + w'$ in the PBW basis, where $\lambda = (\lambda_1 \leq \dots \leq \lambda_d)$ and $w'\in F_{d-1}\calA_c$. By Proposition~\ref{prop: leading term}, $\St_n(D_q)(w')=0$ and 
    \begin{equation}\label{eq: degree reduction proof}
        \St_n (D_q) (w) = C_d \sum_{\lambda} \gamma_\lambda \sum_{a: \lambda_a = \mu} \St_{2d-2}(D_{q'}) (u_{\hat{\lambda}_a}) \tag{1.4.2}.
    \end{equation} 
    Apply the associated graded map to \eqref{eq: degree reduction proof} and note that it commutes with $\St_{2d-2}$. For each $a$ such that $\lambda_a=\mu$, we have $\gr_{d-1}(u_{\hat{\lambda}_a}) = x^\lambda / x_\mu$, so the inner sum of the right-hand side of \eqref{eq: degree reduction proof} becomes $\partial_{x_\mu} x^\lambda$. As $\partial_{x_\mu} \gr(w)$ is a nonzero element of $\mc{X}_{M, d-1}$,  $\St_{2d-2}(D_{q'}) \big( \partial_{x_\mu} \gr(w) \big) \neq 0$ by Lemma~\ref{lem:injective}. Lastly, $C_d \neq 0$ by Lemma~\ref{lem: conseq of choice}(3), and the theorem follows.
\end{proof}

\begin{lemma}\label{lem: center of Bc}
    Let $c \in \kk$ and let $R$ be $\calA_c$ or $\calB_c$ and $w \in R$. If $\ad(e_m)(w) = 0$ for all $m$, then $w \in \kk$. In particular, the center of $R$ is $Z(R) = \kk$.
\end{lemma}
\begin{proof}
    Assume $w \notin \kk$ and $\ad(e_m)(w) = 0$ for all $m$. Let $M \in \NN$ be the maximum of the absolute values of the indices appearing in $f = \gr(w)$ as in Notation~\ref{not: degree reduction}. For $m > 2M$, each monomial of $D_m(f)$ has exactly one letter of index $p+m > M$ and all other indices in $[-M,M]$. Hence it determines the monomial of $f$ and the letter it came from. Since $p-m < 0$ and $r-\frac{m}{2} < 0$, we have $D_m(f) \neq 0$ and $D_m(w) \neq 0$, a contradiction.
\end{proof}

As a result of Theorem~\ref{thm: degree reduction}, simplicity of nonzero central quotient follows readily. 
\simpleVir*
\begin{proof}
    It is easy to see that $\calA_0 = \UU(W)$ is not simple. The augmentation ideal $\UU(W)_+$ generated by elements in $W$ is a proper ideal since $\UU(W)/\UU(W)_+ \cong \kk$.

    Let $V \neq 0$ be a subspace of $\calA_c$ stable under the adjoint action of $W$ and let $w \in V$ be an element of minimal order $d$. If $d \geq 1$, by Theorem~\ref{thm: degree reduction}, $0 \neq \St_{2d+2}(D_q)(w) \in V$ of order $d-1$, contradicting minimality of $d$. Thus $d=0$ and $\kk \subseteq V$. Since every two-sided ideal of $\calA_c$ is $\ad(W)$-stable, $\calA_c$ is simple for all $c \neq 0$.
\end{proof}

\begin{remark}\label{rem: bound is tight}
    The bound $(2k-2)\ell+2$ in Theorem~\ref{thm: tensor power} is sharp for $k=2$, and an example is the intermediate series module for $W$. For $\lambda \in \kk$, let $\mathcal{I}_\lambda$ be the intermediate series module $\{v_m\}_{m \in \ZZ}$ where $e_p \cdot v_m = (m + \lambda p) v_{m+p}$. A direct computation, similar to Lemma~\ref{lem: one block}, shows that $\mathcal{I}_\lambda$ satisfies $\St_4$ and 
    \[\St_3(e_{q_1}, e_{q_2}, e_{q_3}) \cdot v_m=\lambda(1-\lambda)\,\Delta(q)\,v_{m+q_1+q_2+q_3} \neq 0 \qquad \text{for all } m\in\ZZ.\]
    By Proposition~\ref{prop: Rosset split} with $n=2\ell+1$, only ordered partitions $\mc{J}$ of shape $(3, 2, \dots, 2)$ contribute to $\St_{2\ell+1}(e_{q_1}, \dots, e_{q_n})$ on $\Sa^\ell(\mathcal{I}_\lambda)$. Choosing $q_i=2^i$ and $u=v_{K_\lambda} v_{2K_\lambda} \cdots v_{\ell K_\lambda}$ for some choice of $K_\lambda > 2^{n+1}$, we can show that the resulting monomials are all distinct with nonzero coefficients. Thus 
    \[\St_{2\ell+1}(e_{q_1}, \dots, e_{q_n}) \cdot u \neq 0 \Leftrightarrow \lambda \notin\{0,1\}.\]
    Since $W \cong \mathcal{I}_{-1}$ as $W$-modules, the smallest $n$ for which $\Sa^\ell(W)$ and $F_\ell \UU(W)$ satisfy $\St_n$ is $2\ell+2$.
\end{remark}

\section{\texorpdfstring{Central quotients of the Neveu--Schwarz and Ramond \\ enveloping algebras}{Neveu--Schwarz and Ramon}}\label{sec: super}

In this section, we extend the results of Section~\ref{sec: vir} to the enveloping algebra of the $N = 1$ super Virasoro algebras, that is the Neveu--Schwarz and Ramond enveloping algebras. The main idea of the proof is similar to the Virasoro case. In fact, most of the time, we only need to apply the alternating adjoint operator $\St_n(D_q)$, where each $D_{q_i}$ is an even derivation, to reduce the order of a monomial 
\[u = e_{m_1} \dots e_{m_{\ell'}} G_{r_{\ell'+1}} \dots G_{r_\ell} \in \calB_c.
\]
Applying $\St_n(D_q)$ on $G_r$ does not yield any central elements so the odd generators do not introduce any complication at each step. Similar to Theorem~\ref{thm: degree reduction}, for every $w$ of order $d$, we replace a letter $e_{\mu}$ appearing in $\gr(w)$ by applying $\St_{2d+2}(D_q)$, thus reducing the order of $w$. This is sufficient since we can consider $\ad(G_r)(w)$ for some $r$ if the leading term of $\gr(w)$ involves only odd generators. Let 
\begin{equation}\label{def: A'}
    A'_k(r; T) = \prod_{\alpha=1}^k (r-\frac{t_\alpha}{2} + \sum_{\beta > \alpha} t_\beta).
\end{equation}
The following lemmas are a trivial modification of Lemmas ~\ref{lem: term 1 ele} and~\ref{lem: one block} so we omit their proofs.
\begin{lemma}\label{lem: NSR term 1 ele}
    For all $\nu_1, \dots, \nu_k \in \ZZ$ and $r \in \ZZ_\epsilon$, we have 
    \begin{equation*}
    \pushQED{\qed} 
        D_{\nu_1} \dots D_{\nu_k} (G_r) = A'_k(r; \nu) \, G_{r + \wt(\nu)}. \qedhere
    \popQED
    \end{equation*}
\end{lemma} 
\begin{lemma}\label{lem: NSR one block}
    Let $J = \{j_1 < \dots < j_k\} \subseteq [n]$ and $r \in \ZZ_\epsilon$. Then 
    \begin{equation}\label{eq: NSR std ad 1}
        \St_J(D_q)(G_r) = \omega'_J (r) G_{r+\wt(J)}, \qquad \text{where} \qquad \omega'_J (r) = \Alt_k A'_k(r; q_J)
    \end{equation}
    Moreover,
    \begin{enumerate}[label=\textup{(\arabic*)}]
        \item if $k \geq 4$, then $\omega'_J (r) = 0$;
        \item if $J = \{i < i'\}$, then $\omega'_J(r) =(q_{i'} - q_i) (r - \frac{q_i+ q_{i'}}{2})$. \qed
    \end{enumerate} 
\end{lemma}

In other words, the $W$-representation $\mc{G} = \text{span}_{\kk}\{G_r \mid r \in \ZZ_\epsilon\}$ by the adjoint action satisfies the standard identity $\St_4$. The same standard identity is satisfied by the adjoint representation of $W$ over itself, so the odd parts $\mc{G}$ behave exactly like $W$ in this respect. The difference is that a block of four actions on an even letter $e_p$ can leave a nonzero scalar, while it always vanishes on $\mc{G}$.

\begin{proposition}\label{prop: NSR monomial applied}
    Let $c \in \kk$ and $q = (q_1, \dots, q_n) \in \ZZ^n$ with $n = 2h$ even. Let 
    \[u = e_{m_1} \dots e_{m_{\ell'}} G_{r_{\ell'+1}} \dots G_{r_\ell} \in \calB_c.
    \]
    We let $v_a = e_{m_a}$ for $a \leq \ell'$ and $v_a=G_{r_a}$ for $a > \ell'$. Then 
    \begin{equation}\label{eq: NSR monomial applied}
        \St_n(D_q)(u) = \sum_{\mc{J} \in \mc{P}'(u)} \sgn(\mc{J}) \frac{h!}{h_1! \dots h_\ell!} \iprod_{a\in [\ell]} \St_{J_a} (D_q)(v_a),
    \end{equation}
    where $h_a = \length{J_a}/2$ and $\mc{P}'(u)$ is the set of ordered set partitions $\mathcal J=(J_1,\dots,J_\ell)$ such that $\length{J_a} \in \{0, 2, 4\}$ and if $\length{J_a} = 4$, then $a \leq \ell'$ (even part) and $m_a + \wt(J_a) =0$. Moreover, 
    \begin{equation}\label{eq: cases std on v}
        \St_{J_a}(D_q)(v_a) = \begin{cases}
        v_a  & \text{if } J_a = \varnothing; \\
        \omega_{J_a}(m_a)\, e_{m_a + \wt(J_a)} +  \kappa_{J_a}(m_a)\, c,
        & \text{if } \length{J_a} = 2, a \leq \ell'; \\
      \kappa_{J_a}(m_a)\, c, & \text{if } \length{J_a} = 4,a \leq \ell';\\
        \omega'_{J_a}(r_a) \,G_{r_a+\mathrm{wt}(J_a)} & \text{if } \length{J_a} = 2, a > \ell'.
    \end{cases}
    \end{equation}
    In particular, $\St_n(D_q)(u) = 0$ if $n > 2(\ell+ \ell')$ and if $c = 0$ then $\St_n(D_q)(u) = 0$ if $n > 2\ell$.
\end{proposition}
\begin{proof}
    Equations~\eqref{eq: NSR monomial applied} and~\eqref{eq: cases std on v} follow similarly to Proposition~\ref{prop: monomial applied}. We only obtain the new case $\St_{J_a}(D_q)(v_a) =  \omega'_{J_a}(r_a) \,G_{r_a+\mathrm{wt}(J_a)}$ if $\length{J_a} = 2, a > \ell'$. Thus $\St_n(D_q)(u) = 0$ if $n > 4\ell' + 2(\ell-\ell') =  2(\ell+ \ell')$ and if $c = 0$ then $\St_n(D_q)(u) = 0$ if $n > 2\ell$.
\end{proof}

We fix $M$, $n$ and $q$ as in Notation~\ref{not: degree reduction}, where $M$ also bounds $\length{r_a}$. Note that for all $J_a \subseteq \{5, \dots, n\}$ with $\length{J_a} = 2$, and $\length{r} \le M$
\begin{equation}
    \omega'_{J_a}(r) \neq 0,
\end{equation}
since $\frac{1}{2} \wt(J_a) > M$.

\begin{proposition}\label{prop: NSR leading term}
 Let $u = e_{m_1} \dots e_{m_{\ell'}} G_{r_{\ell'+1}} \dots G_{r_\ell}$ be an ordered monomial with all $\length{m_a}$ and  $ \length{r_a} \leq M$, and use Notation~\ref{not: degree reduction}. Then $\St_{2d+2}(D_q)(u) = 0$ if $\ell <d$ and 
    \begin{equation}
            \St_{2d+2}(D_q)(u) = C_d \sum_{a \leq \ell': m_a = \mu} \St_{2d-2}(D_{q'})(u_{\hat{a}}) \qquad \text{if} \qquad \ell=d.
        \end{equation}    
    In particular, $\St_{2d+2}(D_q)(u) \in F_{d-1} \calB_c$ and $\St_{2d+2}(D_q)(u) = 0$ if there exists no $a$ such that $m_a = \mu$.
\end{proposition}
\begin{proof}
    By Proposition~\ref{prop: NSR monomial applied}, if $\length{J_a} = 4$, then it has to be in the even part, i.e., $a \leq \ell'$. Thus similar to Proposition~\ref{prop: leading term}, Lemma~\ref{lem: conseq of choice} implies $J_a = \{1, 2, 3, 4\}$. On the odd part, there is no block $J_a$ of size $4$ and  arguing exactly like Proposition~\ref{prop: leading term} with $\omega'_{J_a}(r_a) \neq 0$, we obtain the proposition.
\end{proof}

Recall that $\Sa(W \oplus \mc{G})$ denotes the supersymmetric algebra of the Lie superalgebra $W \oplus \mc{G}$. Since $D_m$ preserves the PBW filtration, $\St_{2d-2}(D_{q'})$ similarly induces a linear map 
\[\St_{2d-2}(D_{q'}): \Sa^{d-1}(W \oplus \mc{G}) \to \Sa^{d-1}(W \oplus \mc{G}).\]
The following lemma follows almost exactly like Lemma~\ref{lem:injective} as there will be no odd blocks $J_a$ and $J_b$ yielding the same letter. 
\begin{lemma}\label{lem: NSR injective}
    Use Notation~\ref{not: degree reduction} and let $\mc{X}_{M,d-1} \subseteq \Sa^{d-1}(W \oplus \mc{G})$ be the span of monomials of degree $d-1$ with indices bounded by $M$, i.e.,  
    \[ x_{k_1} \dots x_{k_{d'-1}} y_{r_{d'}} \dots y_{r_{d-1}}  \]
    with $k_1 \leq \dots \leq k_{d'-1}$, $r_{d'} < \dots < r_{d-1}$ and $\length{k_i},\length{r_j} \leq M$. 
    The map $\St_{2d-2}(D_{q'})$ is injective on $\mc{X}_{M, d-1}$.
\end{lemma}
\begin{proof}
    We follow the argument of Lemma~\ref{lem:injective}. Let $X= x_{n_1} \dots x_{n_{d'-1}} y_{r_{d'}} \dots y_{r_{d-1}}$. Equation \eqref{eq: injective help 1} becomes
    \begin{equation*}\label{eq: NSR injective help 1}
        \St_{2d-2}(D_{q'})(X) = (d-1)! \sum_{\widehat{\mc{J}}} \sgn(\widehat{\mc{J}}) \prod_{a=1}^{d'-1}\omega_{J_a}(k_a)  \prod_{a=d'}^{d-1}\omega'_{J_a}(r_a) \prod_{a=1}^{d'-1} x_{k_a + \wt(J_a)}\prod_{a=d'}^{d-1} y_{r_a + \wt(J_a)}.
    \end{equation*}
    Applying Lemma~\ref{lem: conseq of choice}(2) to even and odd letters separately, it follows analogously that the images of distinct monomials in $\mc X_{M,d-1}$ have disjoint supports.

    Moreover, since $r_{d'} < \dots < r_{d-1}$, there are no blocks $J_a$ and $J_b$ such that $y_{r_a + \wt(J_a)} = y_{r_b + \wt(J_b)}$. Thus, the lemma follows exactly from the last paragraph of Lemma~\ref{lem:injective}.
\end{proof}

Thus, our order reduction algorithm follows directly.
\begin{theorem}\label{thm: NSR degree reduction}
    Let $c \neq 0$ and let $w \in \calB_c$ of order $d \geq 1$ with $\partial_{x_{\mu}} \gr(w) \neq 0$ for some $\mu \in \ZZ$. Use Notation~\ref{not: degree reduction}. Then we have, 
    \[\gr( \St_{2d+2} (D_q) (w) ) = C_d \St_{2d-2}(D_{q'})(\partial_{x_\mu} \gr(w)) \neq 0.\]
\end{theorem}
\begin{proof}
    We proceed similarly to Theorem~\ref{thm: degree reduction} with Lemma~\ref{lem: NSR injective} and Proposition ~\ref{prop: NSR leading term}. Since $x_\mu$ is even, no sign appears when it is removed.
\end{proof}

 The hypothesis of Theorem \ref{thm: NSR degree reduction} fails exactly when $\gr(w) \in \bigwedge^d G$. In that case, we can easily create an odd letter by applying $\ad(G_s)$ for some $s \in \ZZ_\epsilon$.

\begin{lemma}\label{lem: handle all odds}
    Let $w \in \calB_c$ of order $d \geq 1$ such that $f = \gr (w) \in \bigwedge^d G$ and let $s \in \ZZ_\epsilon$. Then $\ad(G_s)(w)$ has order $d$ and 
    \[\gr(\ad(G_s)(w)) = - 2\sum_{r} x_{r+s} \partial_{y_r} f.\]
    In particular, $\partial_{x_\mu} \gr(\ad(G_s)(w)) \neq 0$ for some $\mu \in \ZZ$.
\end{lemma}
\begin{proof}
    By \eqref{svir-bracket}, we note that $[G_s,G_r]\equiv-2e_{r+s}\pmod{F_0\mathcal{B}_c}$. Thus 
    \[\ad(G_s): \gr \calB_c \to \gr \calB_c, \quad y_r \mapsto -2x_{r+s},\]
    giving the formula. Note that the letters $x_{r+s}$ are all distinct, and since $\partial_{y_r} f \neq 0$ for some $r$, the sum is nonzero.
\end{proof}

In the following theorem, being stable under the adjoint action of $W$ alone is not enough, e.g., $\mc{G} \subseteq F_1 \calB_c$ is a nonzero $\ad(W)$-stable subspace not containing $\kk$. We need the full (super) adjoint action of $\geps$ from Lemma~\ref{lem: handle all odds}.
\begin{theorem}\label{thm: NSR simple 1}
    Let $\g$ be the Neveu--Schwarz or the Ramond algebra with central element $z$ and let $c \neq 0$. Every nonzero subspace of $\UU(\g)/(z-c)$ that is stable under the (super) adjoint action of $\g$ contains $\kk$.
\end{theorem}

\begin{proof}
    We proceed similarly to Theorem~\ref{thm: simple Vir}. Let $V \neq 0$ be a subspace of $\calB_c$ stable under the adjoint action of $\geps$ and let $w \in V$ be an element of minimal order $d$. Assume $d \geq 1$. If $\gr(w) \in \bigwedge^d \mathcal{G}$, we replace $w$ by $\ad(G_s)(w)$, which is nonzero and has the same order by Lemma~\ref{lem: handle all odds}. By Theorem~\ref{thm: NSR degree reduction}, $\St_{2d+2}(D_q)(w) \in V$ has order exactly $d-1$, contradicting minimality of $d$. Thus $d=0$ and $\kk \subseteq V$.
\end{proof}

Unlike the Lie algebra case, for a Lie superalgebra $\g$, a two-sided ideal $I$ of $\UU(\g)$ need not be closed under the super-adjoint action. Let $w = w_0 + w_1$ where $w_0$ and $w_1$ are the even and odd part of $w$ respectively. Then if $w'$ is odd, 
\[[w', w_0+w_1] = w'(w_0+w_1) - w_0 w' +w_1w'\]
involves the homogeneous component $w_0$ and $w_1$, which need not lie in $I$. Nevertheless, the even derivation $\ad(e_0)$ forces such ideal $I$ to be homogeneous with respect to weight. This problem only arises for the Ramond case as for Neveu--Schwarz algebra, the odd weight lies in $\frac{1}{2} + \ZZ$ and the even weight lies in $\ZZ$, which implies that $I$ is closed under the super adjoint action. Nonetheless, we treat both the Neveu--Schwarz and Ramond in the following theorem.

\begin{theorem}\label{thm: NSR simple 2}
    Let $\geps$ be the Neveu--Schwarz or the Ramond algebra. Then $\UU(\geps)/(z-c)$ is simple $\kk$-algebra if and only if $c \neq 0$. 
\end{theorem}
\begin{proof}
    If $c=0$, the augmentation ideal of $\UU(\geps)/(z)$ is a proper nonzero ideal. For the inverse, it suffices to show for the case $\epsilon = 0$ (Ramond), by the above discussion. Let $\theta$ be the parity automorphism induced by 
    \[\theta(e_m) = e_m, \qquad \theta(z) = z, \qquad \theta(G_r) = -G_r, \qquad m, r \in \ZZ.\] 
    Let $I$ be a nonzero proper ideal of $\calB_c$. Then $I + \theta(I)$ and $I \cap \theta(I)$ are homogeneous with respect to parity. By Theorem~\ref{thm: NSR simple 1}, $\calB_c = I \oplus \, \theta(I)$. Let $1 = \chi+\theta(\chi)$ where $\chi \in I$. Then it is easy to check that $\eta = \chi-\theta(\chi)$ is central and $\theta(\eta) = -\eta$, contradicting Lemma~\ref{lem: center of Bc}.
\end{proof}

\section{Applications, Related and Future Work}\label{sec: app}

\subsection{Work of Feigin, Kanel-Belov and Khoroshkin and Etingof's conjecture}\label{subsec: FKBK work}
While we were finishing this version of the paper, Susan Sierra kindly told us about the recent update of \cite{FeiginKanelKhoroshkin}. It gives a different route to proving that $\calA_c$ is simple by proving Etingof's conjecture \cite[Conjecture 4.8]{FeiginKanelKhoroshkin}. Let us first prove \cite[Conjecture 4.8]{FeiginKanelKhoroshkin} and then compare the methods. 

Let $c \in \kk$ and $R$ be $\calA_c$ or $\calB_c$. Let $\mc{K} = \geps/ \kk z = W \oplus \mc{G}$ if $R = \calB_c$ and $\mc{K} = W$ if $R= \calA_c$. For $r \geq 1$, let $E^c_r = F_r R/F_{r-2} R$. Then we have the following short exact sequence of $W$-modules 
\[0 \longrightarrow \Sa^{r-1}(\mc{K}) \longrightarrow E^c_r \stackrel{\pi}{\longrightarrow} \Sa^{r}(\mc{K}) \longrightarrow 0\]
where $\pi$ is the map $F_r R/F_{r-2} R \to F_r R/F_{r-1} R$. For a $\mc{K}$-submodule $\mc{M} \subseteq  \Sa^r(\mc{K})$ write $E^c_r(\mc{M}) := \pi^{-1}(\mc{M})$. Then Etingof's conjecture \cite[Conjecture 4.8]{FeiginKanelKhoroshkin} is a particular instance of the following Theorem~\ref{thm: super Etingof conj} for the module $\Delta^{2N}\Sa^{r}(W)$ and $\mc{G} = 0$.

\begin{theorem}\label{thm: super Etingof conj}
    Let $R= \calA_c$ with $\mc{K} = W$ or $R = \calB_c$ with $\mc{K} = \geps/ \kk z = W \oplus \mc{G}$ for $\epsilon \in \{0, \frac{1}{2}\}$. Let $c \neq 0$ and $r \geq 1$. For every nonzero $\mc{K}$-submodule $\mc{M}$ of $\Sa^r(\mc{K})$, the extension 
    \[0 \longrightarrow \Sa^{r-1}(\mc{K}) \longrightarrow E^c_r (\mc{M}) \stackrel{\pi}{\longrightarrow} \mc{M} \longrightarrow 0\]
    does not split as a sequence of $W$-modules.
\end{theorem}
\begin{proof}
    Assume the contrary and let $s: \mc{M} \to E^c_r (\mc{M})$ be a  section. Let $0 \neq f \in \mc{M}$, and we can assume $f \notin \bigwedge^{r} (\mc{G})$ by possibly replacing $f$ by $G_s \cdot f \in \mc{M}$, which is nonzero by Lemma~\ref{lem: handle all odds}. Let $w \in R$ be a lift of $s(f)$. Then $\gr(w) = \pi(s(f)) = f \neq 0$ so $w$ has order $r \geq 1$. 
    Let $\mu \in \ZZ$ such that $x_\mu$ appears in $f = \gr(w)$ and choose $M, q$ for $w$ as in Notation~\ref{not: degree reduction}. We consider $\St_{2r+2}(D_q) s(f)$.

    On the one hand, the $W$-module $\mc{K}$ satisfies the standard polynomial $\St_4$ by Lemma~\ref{lem: one block}(1) for $W$ and Lemma~\ref{lem: NSR one block}(1) for $\mc{G}$. By Theorem~\ref{thm: tensor power}, the $W$-module $\Sa^r(\mc{K})$ satisfies $\St_{2r+2}$. As $s$ is $W$-equivariant,  
    \[\St_{2r+2}(D_q) \big(s(f)\big) = s(\St_{2r+2}(D_q)(f)) = 0.\]

    On the other hand, by definition, $\St_{2r+2}(D_q)(s(f)) = \St_{2r+2}(D_q)(w)$ as cosets of $F_{r-2} R$. Our order reduction algorithm Theorem~\ref{thm: degree reduction} and Theorem~\ref{thm: NSR degree reduction} imply that  
    \[\gr \big( \St_{2r+2}(D_q)(w) \big) = C_r \, \St_{2r-2}(D_{q'}) \big( \partial_{x_\mu} \gr(w) \big) \neq 0, \]
    so $\St_{2r+2}(D_q)(w)$ has order $r-1$ and $\St_{2r+2}(D_q)(s(f)) \neq 0$ in $E^c_r$, a contradiction.
\end{proof}

\begin{remark}
In \cite{FeiginKanelKhoroshkin}, the link from Etingof's Conjecture \cite[Conjecture 4.8]{FeiginKanelKhoroshkin} to Theorem~\ref{thm: simple Vir} is through the existence of a submodule $\Delta^{2N}\Sa^r(W)$ in every nonzero $W$-module of $\Sa^r(W)$ \cite[Theorem~3.8]{FeiginKanelKhoroshkin}. The main similarity between our paper and \cite{FeiginKanelKhoroshkin} is that both use some versions of alternating combinations of the basis $e_m$ of $W$, but our focus and techniques differ significantly. In \cite{FeiginKanelKhoroshkin}, the operators are applied to $\Sa^r(\mathcal{I}_\lambda)$, where $\mathcal{I}_\lambda$ is a module of tensor fields. They focused on the submodule structure of $\Sa^r(\mathcal{I}_\lambda)$ and proved various interesting results including Noetherianity and Gelfand's conjecture on the finite-dimensionality of the cohomology of graded subalgebras of finite codimension in Lie algebras of vector fields. However, these modules cannot detect the central charge $c$, so the cocycle has to be treated separately via Etingof's conjecture.

By contrast, we construct an explicit operator $\St_{2r+2}(D_q)$ that acts on $R$ itself. It annihilates $\Sa^r(\mc{K})$ and $\Sa^{r-1}(\mc{K})$ (Theorem~\ref{thm: tensor power}), keeping only the central term, which gives a nonzero leading coefficient and lowers the order by exactly one. These results allow us to prove Etingof's conjecture in its strongest form, for every nonzero submodule (Theorem~\ref{thm: super Etingof conj}), and apply verbatim to the Neveu--Schwarz and Ramond algebras. Our methods therefore do not depend on the structure of the module and can be generalized to other Lie algebras in Section~\ref{subsec: Future work}. On the other hand, our paper does not recover their results on $W$-modules, including Noetherianity of tensor products of density modules and the fact that every nonzero submodule of $\Sa^r(W)$ contains some $\Delta^{2N}\Sa^r(W)$ \cite[Theorem~3.8]{FeiginKanelKhoroshkin}. To the best of our knowledge, these do not follow from polynomial identities alone.
\end{remark}

\subsection{Consequences on the spectra of the enveloping algebra of \texorpdfstring{$\Vir$}{Vir}}\label{subsec: conseq}
Besides the applications to the representation theory of $\Vir$ in Corollary~\ref{cor: ann}, Theorem~\ref{thm: simple Vir} also sheds light on the ideal structure of $\UU(\Vir)$. We have the analog of \cite[Theorem 1.1(3)]{BiswalSierra25} that every nonzero ideal of $\UU(\Vir)$ contains a nonzero element of $\kk[z]$. Hence, every nonzero prime ideal of $\UU(\Vir)$ is either $(z-c)$ for some $c \neq 0$ or the preimage of a prime ideal of $\UU(W) \cong \UU(\Vir)/(z)$.

\begin{theorem}\label{thm: ideal UVir}
    Let $I$ be a nonzero two-sided ideal of $\UU(\Vir)$. Then $I \cap \kk[z] \neq 0$. Moreover, $I$ can be written uniquely as
    \[I = g I_0,\]
    where $g \in \kk[z]$ is monic with $g(0) \neq 0$, and $I_0$ is the preimage of an ideal of $\UU(\Vir)/(z^k)$.
    
    Furthermore if $\kk$ is an algebraically closed field of characteristic zero, then
    \[ \Spec \UU(\Vir) = \{(0)\} \sqcup \{(z-a) \mid a \in \kk^\times\} \sqcup \Spec \UU(W). \]
    If $\kk$ is an uncountable and algebraically closed field of characteristic zero, then
    \[ \Prim \UU(\Vir) = \{(z-a) \mid a \in \kk^\times\} \sqcup \Prim \UU(W). \]
\end{theorem}
\begin{proof}
    By Corollary~\ref{cor: simplicity Vir tensor}, $\UU(\Vir) \otimes_{\kk[z]} \kk(z)$ is simple. Let $S = \kk[z] \setminus \{0\}$ and suppose $I \cap S = \emptyset$. As $\UU(\Vir)$ is free over $\kk[z]$, we have $S^{-1}I \neq 0$, so $S^{-1}I  = S^{-1}\UU(\Vir)$. Thus, there exists $x \in I, s \in S$ such that $1 = s^{-1}x$, so $s = x \in I$, a contradiction. Let $I \cap \kk[z] = (f)$ and write $f = z^k g$ with $g(0) \neq 0$. 

    Let $p$ be a monic irreducible factor of $g$. We will prove by contradiction that $I \subseteq p\UU(\Vir)$. Suppose not, let $\kk' = \kk[z]/(p)$ and let $\lambda \in \kk'$ be the image of $z$, which is nonzero. Since  
    \[\UU(\Vir)/(p) \cong \UU(\Vir_{\kk'})/(z -\lambda),\]
    as $\kk'$-algebras, where $\Vir_{\kk'} = \Vir \otimes_{\kk} \, \kk'$ with central element $z$, it is a simple ring by Theorem~\ref{thm: simple Vir}. So $I + p \UU(\Vir) = \UU(\Vir)$ and  $1 = u + py$ for some $u \in I, y \in  \UU(\Vir)$. Let $\tilde{f} = f/p$, then $\tilde{f} = \tilde{f}u + fy \in I$, which contradicts $I \cap \kk[z] = (f)$.

    So $I = pI'$ where $I' = \{u \mid pu \in I\}$ and $I' \cap \kk[z] = (\tilde{f})$. Repeating the process for $I'$ yields $I = gI_0$ where $I_0 \cap \kk[z] = (z^k)$. If $I = g'I_0'$ is any other such decomposition, then $I \cap \kk[z] = (g'z^i)$ for some $i$, showing that $g' = g$ and thus $I_0' = I_0$.

    If $P$ is a nonzero prime ideal of $\UU(\Vir)$, then $P  \cap \kk[z]$ is a nonzero prime ideal of $\kk[z]$, so it equals $(z - a)$. If $a \neq 0$, Theorem \ref{thm: simple Vir}  gives $P = (z - a)$. If $a = 0$, then $P$ corresponds to a prime of $\UU(W)$. If $\kk$ is uncountable, then by \cite[9.1.8(1)]{McConnellRobson01}, $\UU(\Vir)$ satisfies the Nullstellensatz over $\kk$ and thus $z$ acts by a scalar on every simple module, so $0$ is not primitive. On the other hand, each $(z - a)$ is maximal, hence primitive.
\end{proof}

An important consequence is that the conjecture in \cite[Conjecture 8.6]{pham2025orbit} that the Dixmier map 
\[\widetilde{\Dx}^{\Vir}: \Pprim \Sa(\Vir) \to \Prim \UU(\Vir)\]
is surjective is equivalent to the conjecture that the Dixmier map $\widetilde{\Dx}^{W}: \Pprim \Sa(W) \to \Prim \UU(W)$ is surjective. Here $\Pprim \Sa(\g)$ denotes the Poisson primitive spectrum of $\Sa(\g)$.

\subsection{Future work}\label{subsec: Future work}
Our method for this paper involves only two properties of $\Vir$. First, the action of $W$ satisfies the standard identity $\St_4$ and this order is indeed optimal. Second, the Virasoro cocycle breaks this identity by a nonzero scalar. Theorem~\ref{thm: tensor power} carries these two features to every PBW piece of $\mc{A}_c$ and $\mc{B}_c$, which allows us to construct the operator $\St_{2\ell+2}(D_q)$ that lowers the order of any element by exactly one. We thus expect our method to generalize to Lie (super)algebras satisfying a polynomial identity which fails in its central extension by a scalar. In this subsection, we list several conjectures for the Lie (super)algebras that are important in mathematics and physics.  

The natural next cases are other Virasoro-related (super)algebras. The \emph{twisted Heisenberg-Virasoro algebra} $\mathcal{H}$ is spanned by $\Vir= \Span \{e_m,z_L\}_{m \in \ZZ}$ (with $z_L$ in place of $z$) together with elements $I_m$ ($m\in\ZZ$), $z_I$ and $z_{LI}$, where $z_L,z_I,z_{LI},I_0$ are central and
\[
[e_m,I_n] = nI_{m+n} + (m^2+m)\,\delta_{m+n,0} \, z_{LI}, \qquad 
[I_m,I_n] = m\,\delta_{m+n,0} \, z_I.
\]
As the even part of the $N=2$ superconformal algebra \cite[Example 5.9 c]{Kac98vertex} is a copy of $\mc H$ with $z_L=c$, $z_I=\frac{c}{3}$ and $z_{LI}=0$, the new difficulty in both cases comes from the Heisenberg part. 
\begin{conjecture}[Heisenberg-Virasoro and $N=2$ superconformal algebra]\label{conj: HV-N2}
\leavevmode
\begin{enumerate}
\item If $c_I \neq 0$, then
$\UU(\mc H)/(z_L-c_L,\,z_I-c_I,\,z_{LI}-c_{LI},\,I_0-a)$ is simple for all
$c_L,c_{LI},a\in\kk$.
\item Let $\epsilon \in \{0, \frac{1}{2}\}$ and $\g^{N=2}_\epsilon$ be the $N=2$ superconformal algebra. Then $\UU(\g^{N=2}_\epsilon)/(z_L-c_L)$ is simple if and only if $c_L \neq 0$.
\end{enumerate}
\end{conjecture}

The \emph{$\mathfrak{bms}_3$ algebra} appears as the symmetry algebra of three-dimensional asymptotically flat gravity and is the central extension of $W(0,-1) = W \ltimes \mc{I}_{-1}$, where $\mc{I}_{-1}\cong W$ is the adjoint representation. As before, it contains $\Vir$ (with $z_L$ in place of $z$), together with $M_m$ ($m\in\ZZ$) and central element $z_M$ satisfying
\[
[e_m,M_n]=(n-m)\,M_{m+n}+\delta_{m+n,0}\,\tfrac{m^3-m}{12}\,z_M,\qquad
[M_m,M_n]=0.
\]
The \emph{boundary Carrollian conformal algebra} (BCCA) \cite{BHPTVZ26} is the subalgebra 
\[
\hat{\mf{b}} = \Span \{e_n-e_{-n}, \, M_m+M_{-m}, \,  z_M \mid n \geq 1,\ m \geq 0\}
\]
of $\mf{bms}_3$, the algebra of symmetries on the worldsheet of the tensionless open strings.
\begin{conjecture}[$\mf{bms}_3$ and Carrollian algebras]\label{conj: BMS-BCCA}
\leavevmode
\begin{enumerate}
\item The algebra $\UU(\mf{bms}_3)/(z_L-c_L,\,z_M-c_M)$ is simple if and only if $c_M \neq 0$.
\item The algebra $\UU(\hat{\mf{b}})/(z_M-c_M)$ is simple if and only if $c_M \neq 0$.
\end{enumerate}
\end{conjecture}
We also expect the results to hold for the affine--Virasoro algebra $\Vir \ltimes \hat{\g}$ by combining our results and \cite{BiswalSierra25}.

On the other hand, Razmyslov showed that apart from $\mf{sl}_2$, all simple Lie algebras satisfying the standard Lie identity of degree $5$ (equivalently, satisfies $\St_4$ as an adjoint module over itself) are Lie algebras of vector fields over commutative differential algebras \cite{Razmyslov85}. A natural next step is therefore vector fields on curves. 

Let $C$ be a smooth irreducible affine curve with coordinate ring $\mc{O}(C)$ and let $L(C)=\Der \mc{O}(C)$. Each $\varphi \in H^1_{\mathrm{dR}}(C)^*$ corresponds to a Virasoro-type cocycle, given locally by $c_\varphi(f\partial_t,g\partial_t)=\varphi(f'''g\,dt)$, and hence a
one-dimensional central extension $L_\varphi(C)= L(C) \oplus \kk z$
\cite[Chapter 6]{Schlichenmaier14}.

\begin{conjecture}[Vector fields on curves]\label{conj: curves}
Let $C$ be a smooth irreducible affine curve and $\varphi\in H^1_{\mathrm{dR}}(C)^*$. Let $L_\varphi(C) = L(C) \oplus \kk z_L$ be the central extension defined by the cocycle $c_\varphi$. The following are equivalent:
\begin{enumerate}
    \item the class $\varphi \in H^1_{\mathrm{dR}}(C)^*$ is nonzero;
    \item the central quotient $\Sa(L_\varphi(C))/(z_L-1)$ is Poisson simple;
    \item the central quotient $\UU(L_\varphi(C))/(z_L-1)$ is simple.
\end{enumerate}
\end{conjecture}

For $C= \mathbb{G}_m$ and $\varphi \neq 0$, we recover the Virasoro algebra. Then $(1) \Leftrightarrow (2)$ is \cite[Corollary 3.3.5]{PetukhovSierra2022} and $(1) \Leftrightarrow (3)$ is Theorem~\ref{thm: simple Vir}.

\bibliographystyle{amsalpha}
\bibliography{refs.bib}

@article {AmitsurLevitzki50,
    AUTHOR = {Amitsur, S. A. and Levitzki, J.},
     TITLE = {Minimal identities for algebras},
   JOURNAL = {Proc. Amer. Math. Soc.},
  FJOURNAL = {Proceedings of the American Mathematical Society},
    VOLUME = {1},
      YEAR = {1950},
     PAGES = {449--463},
}

@article {Rosset76,
    AUTHOR = {Rosset, Shmuel},
     TITLE = {A new proof of the {A}mitsur-{L}evitski identity},
   JOURNAL = {Israel J. Math.},
  FJOURNAL = {Israel Journal of Mathematics},
    VOLUME = {23},
      YEAR = {1976},
    NUMBER = {2},
     PAGES = {187--188},
    }

@article {BiswalSierra25,
    AUTHOR = {Biswal, Rekha and Sierra, Susan J.},
     TITLE = {Ideals in enveloping algebras of affine {K}ac-{M}oody
              algebras},
   JOURNAL = {Algebra Number Theory},
  FJOURNAL = {Algebra \& Number Theory},
    VOLUME = {19},
      YEAR = {2025},
    NUMBER = {6},
     PAGES = {1199--1230},
     }

@misc{pham2025orbit,
  title         = {The orbit method for the {V}irasoro algebra},
  author        = {Pham, Tuan Anh},
  year          = {2025},
  eprint        = {2504.14670},
  archiveprefix = {arXiv},
  primaryclass  = {math.RA},
  note          = {\href{https://arxiv.org/abs/2504.14670}{arXiv:2504.14670}}
}

@misc{FeiginKanelKhoroshkin,
      title={On finite dimensionality of homology of subalgebras of vector fields}, 
      author={Boris Feigin and Alexei Kanel-Belov and Anton Khoroshkin},
      year={2026},
      eprint={2211.08510},
      archivePrefix={arXiv},
      primaryClass={math.QA},
      note          = {\href{https://arxiv.org/abs/2211.08510}{arXiv:2211.08510v3}}
}

@article {Kostant81,
    AUTHOR = {Kostant, Bertram},
     TITLE = {A {L}ie algebra generalization of the {A}mitsur-{L}evitski
              theorem},
   JOURNAL = {Adv. in Math.},
  FJOURNAL = {Advances in Mathematics},
    VOLUME = {40},
      YEAR = {1981},
    NUMBER = {2},
     PAGES = {155--175},
}

@article{PetukhovSierra2022,
  author   = {Petukhov, Alexey V. and Sierra, Susan J. },
  title    = {The {P}oisson spectrum of the symmetric algebra of the
              {V}irasoro algebra},
  journal  = {Compos. Math.},
  fjournal = {Compositio Mathematica},
  volume   = {159},
  year     = {2023},
  number   = {5},
  pages    = {933--984}
}

@article {Razmyslov85,
    AUTHOR = {Razmyslov, Yu. P.},
     TITLE = {Simple {L}ie algebras that satisfy the standard {L}ie identity
              of degree {$5$}},
   JOURNAL = {Izv. Akad. Nauk SSSR Ser. Mat.},
  FJOURNAL = {Izvestiya Akademii Nauk SSSR. Seriya Matematicheskaya},
    VOLUME = {49},
      YEAR = {1985},
    NUMBER = {3},
     PAGES = {592--634},
}

@article{IyuduSierra20,
  author   = {Iyudu, Natalia K. and Sierra, Susan J.},
  title    = {Enveloping algebras with just infinite {G}elfand-{K}irillov
              dimension},
  journal  = {Ark. Mat.},
  fjournal = {Arkiv f\"or Matematik},
  volume   = {58},
  year     = {2020},
  number   = {2},
  pages    = {285--306}
}

@article {MazorchukZhao14,
    AUTHOR = {Mazorchuk, Volodymyr and Zhao, Kaiming},
     TITLE = {Simple {V}irasoro modules which are locally finite over a
              positive part},
   JOURNAL = {Selecta Math. (N.S.)},
  FJOURNAL = {Selecta Mathematica. New Series},
    VOLUME = {20},
      YEAR = {2014},
    NUMBER = {3},
     PAGES = {839--854},
}

@article {ChenGuoZhao13,
    AUTHOR = {Chen, Hongjia and Guo, Xiangqian and Zhao, Kaiming},
     TITLE = {Tensor product weight modules over the {V}irasoro algebra},
   JOURNAL = {J. Lond. Math. Soc. (2)},
  FJOURNAL = {Journal of the London Mathematical Society. Second Series},
    VOLUME = {88},
      YEAR = {2013},
    NUMBER = {3},
     PAGES = {829--844},
     }

@article {TanZhao16,
    AUTHOR = {Tan, Haijun and Zhao, Kaiming},
     TITLE = {Irreducible {V}irasoro modules from tensor products},
   JOURNAL = {Ark. Mat.},
  FJOURNAL = {Arkiv f\"or Matematik},
    VOLUME = {54},
      YEAR = {2016},
    NUMBER = {1},
     PAGES = {181--200},
     }

@article {KytolaRidout09,
    AUTHOR = {Kyt\"ol\"a, Kalle and Ridout, David},
     TITLE = {On staggered indecomposable {V}irasoro modules},
   JOURNAL = {J. Math. Phys.},
  FJOURNAL = {Journal of Mathematical Physics},
    VOLUME = {50},
      YEAR = {2009},
    NUMBER = {12},
     PAGES = {123503, 51},
     }

@article {CJORY21,
    AUTHOR = {Creutzig, Thomas and Jiang, Cuipo and Orosz Hunziker,
              Florencia and Ridout, David and Yang, Jinwei},
     TITLE = {Tensor categories arising from the {V}irasoro algebra},
   JOURNAL = {Adv. Math.},
  FJOURNAL = {Advances in Mathematics},
    VOLUME = {380},
      YEAR = {2021},
     PAGES = {Paper No. 107601, 35},
}

@book {Kac98vertex,
    AUTHOR = {Kac, Victor},
     TITLE = {Vertex algebras for beginners},
    SERIES = {University Lecture Series},
    VOLUME = {10},
   EDITION = {Second},
 PUBLISHER = {American Mathematical Society, Providence, RI},
      YEAR = {1998},
     PAGES = {vi+201},
    }

@article {BHPTVZ26,
    AUTHOR = {Buzaglo, Lucas and He, Xiao and Pham, Tuan Anh and Tan, Haijun
              and Vishwa, Girish S. and Zhao, Kaiming},
     TITLE = {On the boundary {C}arrollian conformal algebra},
   JOURNAL = {Lett. Math. Phys.},
  FJOURNAL = {Letters in Mathematical Physics},
    VOLUME = {116},
      YEAR = {2026},
    NUMBER = {4},
     PAGES = {Paper No. 86, 62},
     }

@book {Schlichenmaier14,
    AUTHOR = {Schlichenmaier, Martin},
     TITLE = {Krichever-{N}ovikov type algebras},
    SERIES = {De Gruyter Studies in Mathematics},
    VOLUME = {53},
      NOTE = {Theory and applications},
 PUBLISHER = {De Gruyter, Berlin},
      YEAR = {2014},
     PAGES = {xvi+360},
     }

@book{McConnellRobson01,
  author    = {McConnell, John C. and Robson, J. Chris},
  title     = {Noncommutative {N}oetherian rings},
  series    = {Graduate Studies in Mathematics},
  volume    = {30},
  edition   = {Revised},
  publisher = {American Mathematical Society, Providence, RI},
  year      = {2001},
  pages     = {xx+636}
}

@article{SierraWalton14,
  author   = {Sierra, Susan J. and Walton, Chelsea},
  title    = {The universal enveloping algebra of the {W}itt algebra is not
              noetherian},
  journal  = {Adv. Math.},
  fjournal = {Advances in Mathematics},
  volume   = {262},
  year     = {2014},
  pages    = {239--260}
}

@article{SierraWalton16,
  author   = {Sierra, Susan J. and Walton, Chelsea},
  title    = {Maps from the enveloping algebra of the positive {W}itt
              algebra to regular algebras},
  journal  = {Pacific J. Math.},
  fjournal = {Pacific Journal of Mathematics},
  volume   = {284},
  year     = {2016},
  number   = {2},
  pages    = {475--509}
}

@misc{PhamTimmins25,
  title         = {The Kernel and Image of Orbit Homomorphisms for the {W}itt Algebra},
  author        = {Pham, Tuan Anh and Timmins, James},
  year          = {2025},
  eprint        = {2510.00756},
  archiveprefix = {arXiv},
  primaryclass  = {math.RA},
  note          = {\href{https://arxiv.org/abs/2510.00756}{arXiv:2510.00756}}
}

@article {Latyvsev63,
    AUTHOR = {Laty{\v s}ev, V. N.},
     TITLE = {Two remarks on {PI}-algebras},
   JOURNAL = {Sibirsk. Mat. \v Z.},
  FJOURNAL = {Akademija Nauk SSSR. Sibirskoe Otdelenie. Sibirski\u i\
              Matemati\v ceski\u i\ \v Zurnal},
    VOLUME = {4},
      YEAR = {1963},
     PAGES = {1120--1121},
      ISSN = {0037-4474},
      }

@article {Bahturin74,
    AUTHOR = {Bahturin, Yuri},
     TITLE = {Identities in the universal envelopes of {L}ie algebras},
   JOURNAL = {J. Austral. Math. Soc.},
  FJOURNAL = {Journal of the Australian Mathematical Society},
    VOLUME = {18},
      YEAR = {1974},
     PAGES = {10--21},}

@article {SanchezSierra23,
    AUTHOR = {Le\'on S\'anchez, Omar and Sierra, Susan J.},
     TITLE = {A {P}oisson basis theorem for symmetric algebras of
              infinite-dimensional {L}ie algebras},
   JOURNAL = {Ark. Mat.},
  FJOURNAL = {Arkiv f\"or Matematik},
    VOLUME = {61},
      YEAR = {2023},
    NUMBER = {2},
     PAGES = {375--412},
}

@misc{pham2026ideals,
      title={Ideals of homomorphic images of the enveloping algebra of the {W}itt algebra}, 
      author={Pham, Tuan Anh},
      year={2026},
      eprint={2610.01913},
      archivePrefix={arXiv},
      primaryClass={math.RA},
      note          = {\href{https://arxiv.org/abs/2610.01913}{arXiv:2610.01913}}
}

@book{Razmyslov94,
    AUTHOR = {Razmyslov, Yu. P.},
     TITLE = {Identities of algebras and their representations},
    SERIES = {Translations of Mathematical Monographs},
    VOLUME = {138},
      NOTE = {Translated from the 1989 Russian original by A. M. Shtern},
 PUBLISHER = {American Mathematical Society, Providence, RI},
      YEAR = {1994},
     PAGES = {xiv+318},
     }

@article {Drensky21,
    AUTHOR = {Drensky, Vesselin},
     TITLE = {Weak polynomial identities and their applications},
   JOURNAL = {Commun. Math.},
  FJOURNAL = {Communications in Mathematics},
    VOLUME = {29},
      YEAR = {2021},
    NUMBER = {2},
     PAGES = {291--324},
     }

@article {Regev72,
    AUTHOR = {Regev, Amitai},
     TITLE = {Existence of identities in {$A\otimes B$}},
   JOURNAL = {Israel J. Math.},
  FJOURNAL = {Israel Journal of Mathematics},
    VOLUME = {11},
      YEAR = {1972},
     PAGES = {131--152},
      }
\end{document}